\pdfoutput=1
\documentclass{jams-l}
\usepackage[francais]{babel}
\usepackage[latin1]{inputenc}
\usepackage{amssymb}
\usepackage[mathscr]{euscript}
\usepackage{rotating}
\usepackage[all]{xy}
\usepackage{hyperref}

\newtheorem{theoreme}[subsubsection]{Théorème}
\newtheorem{lemme}[subsubsection]{Lemme}
\newtheorem{proposition}[subsubsection]{Proposition}
\newtheorem{corollaire}[subsubsection]{Corollaire}
\newtheorem{utheoreme}{Théorème}

\theoremstyle{definition}
\newtheorem{definition}[subsubsection]{Définition}
\newtheorem{notation}[subsubsection]{Notation}

\theoremstyle{remark}
\newtheorem{remarques}[subsubsection]{Remarques}
\newtheorem{remarque}[subsubsection]{Remarque}

\newcommand\C{\mathbb{C}}
\newcommand\Nat{\mathbb{N}}
\newcommand\R{\mathbb{R}}
\newcommand\Z{\mathbb{Z}}
\newcommand\Q{\mathbb {Q}}
\newcommand\Fi{\mathbb{F}}
\newcommand\Af{\mathbb{A}_f}

\newcommand\Hp{{}^{p}{H}}

\DeclareMathOperator\Hom{\mathrm{Hom}}
\DeclareMathOperator\Ext{\mathrm{Ext}}

\newcommand\SD{\mathbb{S}}

\newcommand\G{{\bf G}}
\newcommand\GU{{\bf GU}}
\newcommand\GSp{{\bf GSp}}
\newcommand\Sp{{\bf Sp}}
\newcommand\GL{{\bf GL}}
\newcommand\Pa{{\bf P}}
\newcommand\QP{{\bf Q}}
\newcommand\RP{{\bf R}}
\newcommand\B{{\bf B}}
\newcommand\Gr{\mathbb{G}}
\newcommand\N{{\bf N}}
\newcommand\U{{\bf U}}
\newcommand\Le{{\bf L}}
\newcommand\Se{{\bf S}}

\newcommand\K{\mathrm {K}}
\newcommand\Hr{\mathrm{H}}

\newcommand\F{{\mathcal F}}

\DeclareMathOperator{\Proj}{Proj}
\DeclareMathOperator{\card}{card}
\newcommand\DP{{}^{w}D}
\newcommand\as{\underline{a}}
\newcommand\lin{\ell}
\newcommand\Sc{{\mathcal S}}
\newcommand\Mod{{\mathcal M}}
\newcommand\sous{\setminus}
\newcommand\til{\widetilde}
\newcommand\X{{\mathcal X}}
\newcommand\Y{{\mathcal Y}}
\newcommand\fl{\rightarrow}
\newcommand\fle{\mapsto}
\newcommand\iso{\stackrel {\sim} {\fl}}
\newcommand\limid {\displaystyle{\lim_{\overrightarrow{\scriptstyle \,\,\,\Delta\,\,\,}}}}

\title[Complexes d'intersection des compactifications de Baily-Bore]
{Complexes d'intersection des compactifications de Baily-Borel :
Le cas des vari\'et\'es de Siegel}
\author{Sophie Morel}
\date{}
\address{Laboratoire de math\'ematique, Universit\'e Paris Sud, b\^atiment 425, 91405 Orsay Cedex, France}
\curraddr{Princeton Mathematics Department, Fine Hall, Washington Road, Princeton NJ 08544-1000, USA}
\email{smorel@math.princeton.edu}

\begin{document}

\maketitle

\tableofcontents

\begin{abstract}
In this work, we calculate the trace of a Hecke correspondance
composed with a power of the Frobenius endomorphism on the fibre of the
intersection complexes of the Baily-Borel compactification of a Siegel modular variety.

Our main tool is Pink's theorem about the restriction to the strata of the Baily-Borel
compactification of the direct image of a local system on the Shimura variety.
To use this theorem, we give a new construction of the intermediate extension
of a pure perverse sheaf as a weight truncation of the full direct image.

More generally, we are able to define analogs in positive characteristic of the
weighted cohomology complexes introduced by Goresky, Harder and MacPherson.

\end{abstract}

Les variétés de Shimura les plus étudiées sont celles associées au groupe $\GL_2$,
autrement dit les courbes modulaires.
Si $Y$ est une courbe modulaire, on obtient sa compactification de Baily-Borel $j:Y\fl Y^*$
en ajoutant un nombre fini de pointes
(car $\GL_2$ est de rang semi-simple $1$).
Comme $Y^*$ est lisse, le complexe d'intersection associé à un système de coefficients
$\F$ sur $Y$ est $j_*\F$.
Il est possible dans ce cas de calculer la fonction $L$ de $Y^*$ à coefficients dans $j_*\F$,
et il a été prouvé dans des travaux d'Eichler, Shimura, Deligne et Ihara (entre autres)
qu'elle s'écrit comme un produit alterné de fonctions $L$ de formes modulaires cuspidales
et de fonctions zêta.

La fonction $L$ est un produit de facteurs locaux $L_p$, où $p$ parcourt l'ensemble
des nombres premiers, et ce sont les $L_p$ que l'on calcule.
Le cas essentiel est celui où $Y^*$ et $j_*\F$ ont bonne réduction en $p$. Le facteur
local $L_p$ ne dépend alors que des réductions modulo $p$ de $Y^*$ et $j_*\F$.

Pour calculer $L_p$ dans le cas des courbes modulaires, il existe deux méthodes :
la méthode des congruences, qui ne se généralise pas en dimension supérieure,
et la comparaison de la formule des traces d'Arthur-Selberg 
et de la formule des points fixes de Grothendieck-Lefschetz.
C'est cette deuxième méthode que l'on cherche à généraliser.

Pour une variété de Shimura générale $M^{\K}(\G,\X)$ , l'application de cette méthode 
au calcul de $L_p$ est plus délicate.
Un premier pas 
est le calcul de la trace d'une correspondance de Hecke composée avec une puissance du morphisme de Frobenius
sur la cohomologie du complexe d'intersection 
de la compactification de Baily-Borel, ou,
ce qui suffit 
grâce à la conjecture de Deligne si la
puissance de Frobenius est assez grande (cf [P3], [F] et [V]),
des termes locaux naïfs de ces correspondances.

Brylinski et Labesse ([BL]) ont effectué ce calcul 
pour $\G=R_{E/\Q}\GL_2$, avec $E$ une extension totalement réelle 
de $\Q$ de degré supérieur ou égal à $2$, et ils en ont déduit
que la fonction $L$ du complexe d'intersection
était bien de la forme attendue.

Pour $\G=\GU(2,1)$, le calcul a été fait par Kottwitz et Rapoport dans l'article [KR] du livre [LR] (dans ce cas, on peut aussi montrer que la fonction $L$ à coefficients dans
le complexe d'intersection est un produit alterné de fonctions $L$ automorphes,
voir l'article de Langlands et Ramakrishnan dans le même livre [LR]).
Le cas d'un groupe $\G$ de rang $1$ et d'une correspondance de Hecke triviale a été traité par Rapoport dans son article [R], paru aussi dans [LR].

Signalons enfin que le cas d'un groupe unitaire sur $\Q$ de rang quelconque et d'une
correspondance de Hecke triviale a été traité dans [M].

Dans cet article, nous calculons les termes locaux naïfs d'une correspondance de Hecke
composée avec une puissance du morphisme de Frobenius 
pour le complexe d'intersection de la compactification de Baily-Borel des variétés
de Shimura associées aux groupes symplectiques sur $\Q$ (de rang arbitraire).

\vspace{1cm}

Notre outil principal est le théorème de Pink calculant les restrictions aux strates
de la compactification de Baily-Borel du prolongement d'un système de coefficients 
sur la variété de Shimura (cf [P2]).
Pour les variétés de Shimura $M^{\K}(\G,\X)$ considérées dans cet article,
le théorème de Pink s'écrit (cf le théorème \ref{th_Pink})
$$i^*Rj_*\F^{\K}V\simeq\F^{\K'}R\Gamma(\Gamma_\lin,R\Gamma(Lie(\N),V)).$$
$j$ est l'inclusion de la variété de Shimura dans sa compactification de Baily-Borel $M^{\K}(\G,\X)^*$,
$V$ est une représentation algébrique du groupe $\G$,
$\F^{\K}V$ est le système de coefficients associé à $V$,
$i$ est l'inclusion d'une strate associée à un parabolique maximal $\Pa$,
$\N$ est le radical unipotent de $\Pa$
et $\Gamma_\lin$ est un sous-groupe arithmétique de la parte linéaire du quotient de Levi de $\Pa$.

\vspace{1cm}

Pour pouvoir utiliser ce théorème, nous donnons une nouvelle construction du
prolongement intermédiaire d'un faisceau pervers pur :
Si $j:U\fl X$ est l'inclusion d'un ouvert non vide dans un schéma $X$ 
séparé de type fini sur un corps fini
et $K$ un faisceau pervers pur de poids $a$ sur $U$, alors
$$j_{!*}K=w_{\leq a}Rj_*K.$$
Dans cette formule, $w_{\leq a}$ est le tronqué pour la t-structure $(\DP^{\leq a}(X),\DP^{>a}(X))$,
où $\DP^{\leq a}(X)$ (resp. $\DP^{>a}(X)$) est la sous-catégorie pleine de la catégorie des complexes mixtes sur $X$
dont les objets sont les complexes qui ont tous leurs faisceaux de cohomologie perverse
de poids $\leq a$ (resp. $>a$).
Cette t-structure est assez inhabituelle, puisqu'elle est dégénérée et de coeur nul.

De plus, la t-structure $(\DP^{\leq a},\DP^{\geq a+1})$ sur un schéma stratifié
s'obtient en recollant les t-structures analogues sur les strates.

\vspace{1cm}

En combinant la formule ci-dessus et le théorème de Pink, on obtient
le théorème principal de cet article (théorème \ref{th:restriction_W_bord}) :

\begin{utheoreme} On note $c$ la dimension de $M^{\K}(\G,\X)$.
Soit $V$ une représentation algébrique de $\G$ sur laquelle
le centre déployé de $\G$ agit trivialement.
Si $i$ est l'inclusion d'une strate de bord de $M^{\K}(\G,\X)^*$, 
on a l'égalité virtuelle :
$$[i^*IC_V]=
\F^{\K_r}\left(\sum_{S}\sum_{i\in I_S}(-1)^{\card(S)-1}\left[R\Gamma\left(\Gamma_{S},R\Gamma(Lie(\N_{S}),g_ig.V)_{\geq t_{r},< t_{s},s\in S\sous\{r\}}\right)\right]\right),$$
où $IC_V$ est le complexe d'intersection sur $M^{\K}(\G,\X)^*$ à coefficients dans $V$,
défini par la formule suivante :
$$IC_V=(j_{!*}(\F^{\K}V[c]))[-c].$$

\end{utheoreme}

Les notations précises sont expliquées dans 4.2.
Disons seulement 
que les $t_i$ sont des entiers qui ne dépendent que des dimensions des différentes strates
(on peut prendre par exemple $t_i$ égal à l'opposé de la codimension de la $i$-ième strate),
que $S$ parcourt un système d'indices des sous-groupes paraboliques standard
dont les strates de bord associées dans la compactification de Borel-Serre réductive
s'envoient sur la strate de la compactification de Baily-Borel considérée,
que, si $\Pa_{S}$ est le sous-groupe parabolique correspondant à $S$,
alors $\N_{S}$ est le radical unipotent de $\Pa_{S}$ 
et $\Gamma_{S}$ est un sous-groupe arithmétique de la partie linéaire du
quotient de Levi $\Pa_{S}/\N_{S}$ de $\Pa_{S}$,
et enfin que $R\Gamma(Lie(\N_{S}),g_ig.V)_{\geq t_{d-r},<t_{d-s},s\in S\sous\{r\}}$ est un tronqué
pour les poids de certains tores centraux de $\Pa_{S}/\N_{S}$.

L'action des correspondances de Hecke sur les complexes de systèmes locaux 
 qui apparaissent dans ce théorème
est donnée par le théorème \ref{th:restr_corr_sigma}.
Grâce à ce théorème, à un résultat de Kottwitz sur le comptage des points à valeurs dans
un corps fini ([K]) et à la conjecture de Deligne évoquée plus haut,
on peut en déduire une formule pour la trace 
sur la cohomologie du complexe d'intersection à coefficients dans $V$
d'une correspondance de Hecke
composée avec une puissance assez grande du morphisme de Frobenius.

\vspace{1cm}

Du point de vue topologique, 
le complexe d'intersection peut être calculé en utilisant la compactification de Borel-Serre réductive.
Plus précisément, Goresky, Harder et MacPherson construisent une famille de complexes pondérés
sur la compactification de Borel-Serre réductive
et montrent que le complexe d'intersection sur la compactification de Baily-Borel est l'image 
directe de deux de ces complexes pondérés
par le morphisme naturel 
de la compactification de Borel-Serre réductive sur la compactification de Baily-Borel (cf [GHM]).

Nous définissons en caractéristique finie des complexes pondérés analogues
aux images directes des complexes pondérés de Goresky, Harder et MacPherson, 
en remplaçant les troncatures par les poids pour les actions
des tores centraux des groupes associés aux composantes de bord
par les troncatures par les poids de l'endomorphisme de Frobenius sur les strates de bord.
Les premiers, à notre connaissance, à utiliser ce type de méthode sont Looijenga et Rapoport
dans [LR2].
Le théorème ci-dessus vaut en fait pour un complexe pondéré quelconque.

\vspace{1cm}

Passons en revue les différentes parties.

Les deux premières parties contiennent des rappels.

Dans la partie 1, 
nous introduisons les variétés de Shimura que nous comptons étudier et leurs modèles entiers,
puis nous rappelons les théorèmes d'existence des modèles canoniques sur $\Q$ de la compactification
de Baily-Borel et des compactifications toroïdales.
En suivant Pink ([P2] 3.7), nous définissons une stratification du bord
de la compactification de Baily-Borel par des variétés de Shimura associées à des groupes symplectiques plus petits.
Enfin, nous rappelons rapidement certains des résultats de Chai et Faltings
sur les modèles entiers des compactifications.

La partie 2 présente la construction des systèmes de coefficients sur la variété de Shimura
provenant de représentations du groupe.
Nous rappelons d'abord la construction de ces systèmes de coefficients sur les points complexes,
puis nous expliquons une méthode, due à Pink ([P2] 1),
pour montrer que ces systèmes de coefficients proviennent de faisceaux étales sur les modèles canoniques.
Ensuite, nous énonçons le théorème de Pink sur le prolongement de ces faisceaux étales 
à la compactification de Baily-Borel.

La partie 3 est indépendante des autres.
Dans un premier temps, nous y étudions la t-structure $(\DP^{\leq a}(X),\DP^{>a}(X))$ définie plus haut
et y montrons la formule pour le prolongement intermédiaire 
d'un faisceau pervers pur $K$ de poids $a$ sur un ouvert non vide $j:U\fl X$.
Dans un deuxième temps, nous considérons des t-structures sur un schéma stratifié $X$ 
qui s'obtiennent en recollant des t-structures $(\DP^{\leq a'},\DP^{>a'})$ sur les strates (avec un $a'$ qui dépend de la strate).
Grâce aux propriétés de ces t-structures,
nous obtenons si l'ouvert $U$ est réunion de strates une égalité
entre les classes dans le groupe de Grothendieck d'un $w_{\leq a}Rj_*K$ 
et d'une somme alternée de tronqués par le poids qui se calculent sur les strates contenues dans $X-U$ (théorème \ref{th:simplification_formule_traces}).

Dans la partie 4, nous appliquons les résultats de la partie 3 aux variétés de Shimura.
Dans 4.1, nous définissons les complexes pondérés
et montrons que deux de ces complexes sont isomorphes au complexe d'intersection.
Dans 4.2, à l'aide du théorème de Pink et du théorème \ref{th:simplification_formule_traces},
nous obtenons une formule explicite, dans le groupe de Grothendieck,
pour la restriction à une strate de bord de la compactification de Baily-Borel
d'un complexe pondéré.

Dans la partie 5, nous traitons le cas des correspondances de Hecke.
Nous commençons par montrer, dans 5.1, quelques résultats généraux sur les correspondances 
cohomologiques entre tronqués par le poids, en particulier un analogue
pour les correspondances du théorème \ref{th:simplification_formule_traces}
(proposition \ref{prop:s_f_t}).
Dans 5.2, nous calculons explicitement les correspondances cohomologiques qui apparaissent
au second membre dans la proposition \ref{prop:s_f_t}.
Ce calcul ressemble beaucoup à celui fait dans 4.2, mais il faut d'abord 
compléter le théorème de Pink par la proposition \ref{prop:suite_Pink}.

\vspace{1cm}

C'est un plaisir de remercier G. Laumon, qui a passé beaucoup de temps à discuter
avec moi.
Je remercie également M. Harris pour ses critiques d'une première version
de mon texte, et
R. Kottwitz, qui a corrigé ou simplifié les démonstrations de certains résultats
de la partie 3.

\section{Variétés de Shimura et leurs compactifications}

\subsection{Variétés de Shimura}
\hspace{.5cm}

Pour tout $d\in\Nat$, on note
$$J_{d}=\left(\begin{array}{ccc}0 & & 1 \\ & \begin{turn}{45}\large\ldots\end{turn}& \\ 1 & & 0\end{array}\right)\in\GL_{d}(\Z)$$
et
$$J_{d,d}=\left(\begin{array}{cc}0 & J_d \\
-J_d & 0\end{array}\right)\in\GL_{2d}(\Z).$$

Le groupe général symplectique $\GSp_{2d}$ est le schéma en groupes sur $\Z$
dont l'ensemble des points à valeurs dans une $\Z$-algèbre $A$ est
$$\GSp_{2d}(A)=\{g\in\GL_{2d}(A)\mbox{ tq }{}^tgJ_{d,d}g=c(g)J_{d,d},c(g)\in A^\times\}.$$
$\G=\GSp_{2d,\Q}$ est un groupe réductif connexe déployé sur $\Q$, de $\Q$-rang semi-simple $d$, $c:\G\fl\Gr_{m,\Q}$ est un caractère de $\G$,
et le groupe dérivé de $\G$ est $\Sp_{2d,\Q}=ker(c)$.

On note $\X_{d}^+$ ou $\X^+$ (resp. $\X_d^-$ ou $\X^-$) l'ensemble des morphismes
$h:\SD\fl\G_{\R}$ (où $\SD=R_{\C/\R}\Gr_m$ est le tore de Serre)
qui induisent une structure de Hodge pure de type $\{(0,1),(1,0)\}$ sur $\Q^{2d}$,
et tels que la forme bilinéaire $\R^{2d}\times\R^{2d}\fl\R$, $(v,w)\fle {}^tvJ_{d,d}h(i)w$
soit symétrique définie positive (resp. négative).
$\G(\R)$ agit transitivement (par conjugaison) sur $\X=\X^+\cup\X^-$,
et le morphisme
$$h_0:z=a+ib\fle\left(\begin{array}{cc}aI_d & -bJ_d \\
bJ_d & aI_d\end{array}\right)$$
est dans $\X^+$, 
donc $\X\simeq\G(\R)/Stab_{\G(\R)}(h_0)$.

\begin{definition} Soient $d,n\in\Nat^*$,
$S$ un schéma sur $\Z[1/n]$ et
$(A\fl S,\lambda)$ un schéma abélien sur $S$ 
principalement polarisé et de dimension relative $d$.
On note $A[n]$ le schéma en groupes fini étale sur $S$ des points de $n$-torsion de $A$.
Une structure de niveau $n$  
sur $A$ est
un isomorphisme (au-dessus de $S$) $\eta:A[n]\iso \underline{(\Z/n\Z)}^{2d}_S$  
qui, à un automorphisme du $S$-schéma en groupes $\underline{\Z/n\Z}_S$ près,
envoie l'accouplement de Weil $A[n]\times A[n]\fl\underline{\Z/n\Z}_S$ associé à $\lambda$
sur la forme
bilinéaire alternée sur $\underline{(\Z/n\Z)}^{2d}$ de matrice $J_{d,d}$.

\end{definition}

\begin{definition} On note $\Mod_{d,n}$ le champ sur $\Z[1/n]$ des schémas abéliens principalement polarisés 
de dimension relative $d$
avec structure de niveau $n$.

\end{definition}

$\Mod_{d,n}$ est au choix de la forme alternée près le champ défini par Chai et Faltings dans [CF] IV.6.1, mais vu comme champ sur $\Z[1/n]$ 
(au lieu de $\Z[1/n,e^{2i\pi/n}]$, cf la remarque [CF] IV.6.12).
C'est un champ de Deligne-Mumford, lisse et de dimension relative $d(d+1)/2$ sur $\Z[1/n]$.
Si $n\geq 3$, $\Mod_{d,n}$ est un espace algébrique,
et même un schéma quasi-projectif sur $\Z[1/n]$.
La fibre générique de $\Mod_{d,n}$ est la variété de Shimura $M^{\K_d(n)}(\G,\X)$
associée à la donnée de Shimura pure $(\G,\X)$ et au sous-groupe ouvert compact
$\K_d(n)=Ker(\GSp_{2d}(\widehat{\Z})\fl\GSp_{2d}(\Z/n\Z))$ de $\G(\Af)$.
En particulier, l'ensemble des points complexes de $\Mod_{d,n}$ est
$$\Mod_{d,n}(\C)=M^{\K_d(n)}(\G,\X)(\C)=\G(\Q)\sous(\X\times\G(\Af)/\K_d(n)).$$
Comme $\G^{der}=Ker(c)$ et $c:\G\fl\Gr_{m,\Q}$ est surjectif,
on a
$$\pi_0(\Mod_{d,n}(\C))\simeq\R^{+\times}\Q^\times\sous\mathbb{A}^\times/c(\K_d(n))=\Q^{+\times}\sous\Af^\times/c(\K_d(n))\simeq (\Z/n\Z)^\times.$$

Si $d=0$, on a $\G=\Gr_{m,\Q}$,
$\X_0=\{\pm 1\}=\pi_0(\R^\times)$ avec l'action évidente de $\G(\R)=\R^\times$, et
$$\Mod_{0,n}=Spec(\Z[1/n,e^{2i\pi/n}]),$$
vu comme un schéma sur $\Z[1/n]$.

Soit $d\in\Nat$. 
Si $n$ divise $m$, on a un morphisme évident $T_1:\Mod_{d,m}\fl\Mod_{d,n|Spec(\Z[1/m])}$
(qui est un morphisme d'oubli d'une partie de la structure de niveau).
Ce morphisme est fini étale galoisien de groupe $\K_d(n)/\K_d(m)$.
De plus, $\GSp_{2d}(\Z/n\Z)$ agit sur $\Mod_{d,n}$ : 
l'image par $g\in\GSp_{2d}(\Z/n\Z)$ d'un schéma abélien
principalement polarisé de dimension relative $d$
avec structure de niveau $n$ $(A,\lambda,\eta)$
est $(A,\lambda,g\circ\eta)$.

Pour $g\in\GSp_{2d}(\widehat{\Z})$ et pour $n$ divisant $m$, 
on note $T_g:\Mod_{d,m}\fl\Mod_{d,n}$ le composé du morphisme
$T_1:\Mod_{d,m}\fl\Mod_{d,n}$ défini ci-dessus et de l'endomorphisme de $\Mod_{d,n}$
induit par l'image de $g$ dans $\GSp_{2d}(\Z/n\Z)$.

Soit $p$ un nombre premier qui ne divise pas $m$. On note $\Af^p$ l'anneau des adèles finies
dont la composante en $p$ est triviale, et $\Z_{(p)}$ le localisé de $\Z$ en $p$. En utilisant
le problème de modules de \cite{K} 5 au lieu de celui décrit plus haut, on peut définir pour
$g\in\G(\Af^p)$ tel que $g^{-1}\K_d(m)g\subset \K_d(n)$
un morphisme fini étale (de degré $[\K_d(n):g^{-1}\K_d(m)g]$)
$T_g:\Mod_{d,m|\Z_{(p)}}\fl\Mod_{d,n|\Z_{(p)}}$ (cf le début de \cite{K}
6), qui est égal à la restriction du morphisme défini ci-dessus si $g\in\GSp_{2d}(\widehat{
\Z})$. Si $g^{-1}\K_d(m)g$ est un sous-groupe distingué de $\K_d(n)$, le morphisme $T_g$
est galoisien de groupe $\K_d(n)/g^{-1}\K_d(m)g$.

Sur les points complexes, ce morphisme $T_g$ est le morphisme
$$\begin{array}{rcl}\G(\Q)\sous (\X\times\G(\Af)/\K_d(m)) & \fl &
\G(\Q)\sous (\X\times\G(\Af)/\K_d(n)) \\
\lbrack(x,h)] & \fle & [(x,hg)].\end{array}$$

Dans la suite, on suppose toujours $n\geq 3$ et on note $\G(\Z)$, $\G(\widehat{\Z})$, etc,
au lieu de $\GSp_{2d}(\Z)$, $\GSp_{2d}(\widehat{\Z})$, etc.

\subsection{Compactifications sur $\Q$}
\hspace{.5cm}

Nous allons énoncer quelques-uns des résultats de Pink ([P1]) dans le cas particulier considéré.
Commençons par rappeler la description des sous-groupes paraboliques maximaux de $\G$.
Le tore diagonal
$$T=\left\{\lambda\left(\begin{array}{ccc}\lambda_1 & 0 & 0 \\
0 & \ddots & 0\\
0 & 0 & \lambda_{2d} \end{array}\right),\lambda,\lambda_1,\dots,\lambda_{2d}\in\Gr_{m,\Q},
\lambda_r\lambda_{2d-r+1}=1\right\}$$
est un tore maximal de $\G$, et l'intersection $\B$ de $\G$ avec le groupe des matrices triangulaires 
supérieures de $\GL_{2d,\Q}$ est un sous-groupe de Borel de $\G$.
On appelle sous-groupes paraboliques standard de $\G$ les sous-groupes paraboliques de $\G$
qui contiennent $\B$.
Tout sous-groupe parabolique de $\G$ est conjugué par $\G(\Q)$ à un unique sous-groupe parabolique standard,
et les sous-groupes paraboliques standard maximaux sont $\Pa_0,\dots,\Pa_{d-1}$, avec
$$\Pa_r=\left\{\left(\begin{array}{ccc}A & * & * \\
0 & B & * \\
0 & 0 & C\end{array}\right),A,C\in\GL_{d-r,\Q},B\in\GSp_{2r,\Q},{}^tAJ_{d-r}C=c(B)J_{d-r}\right\}.$$
Soit $r\in\{0,\dots,d-1\}$. 
On note $\N_r$ le radical unipotent de $\Pa_r$, 
$$\U_r=\left(\begin{array}{ccc}I_{d-r} & 0 & * \\
0 & I_{2r} & 0 \\
0 & 0 & I_{d-r}\end{array}\right)$$
le centre de $\N_r$,
$\Le_r=\Pa_r/\N_r$,
$$\QP_r=\left\{\left(\begin{array}{ccc}c(B)I_{d-r} & * & * \\
0 & B & * \\
0 & 0 & I_{d-r}\end{array}\right),B\in\GSp_{2r,\Q}\right\}\subset\Pa_r$$
si $r>0$,
$$\QP_0=\left(\begin{array}{cc}\Gr_{m,\Q} I_d & * \\
0 & I_d\end{array}\right)\subset\Pa_0,$$
$\G_r=\QP_r/\N_r$ et
$$\Le_{\lin,r}=\left\{\left(\begin{array}{ccc}A & 0 & 0 \\
0 & I_{2r} & 0 \\
0 & 0 & C\end{array}\right),A,C\in\GL_{d-r,\Q},{}^tAJ_{d-r}C=J_{d-r}\right\}\subset\Le_r.$$
On a $\Le_r=\Le_{\lin,r}\times\G_r$, avec
$\Le_{\lin,r}\simeq\GL_{d-r,\Q}$ et $\G_r\simeq\GSp_{2r,\Q}$.
Si $A$ est une $\Z$-algèbre, on note
$\Le_{\lin,r}(A)$ et $\G_r(A)$ les images réciproques par ces isomorphismes
de $\GL_{d-r}(A)$ et $\GSp_{2r}(A)$.

$\QP_r$ est le sous-groupe distingué de $\Pa_r$ défini par Pink dans [P1] 4.7 (cf [P1] 4.25).
On peut donc construire comme dans [P1] 4.11 un espace homogène $\Y_r$ sous $\QP_r(\R)\U_r(\C)$
et une application $h_r:\Y_r\fl Hom(\SD_\C,\QP_{r,\C})$
tels que $(\QP_r,\Y_r)$ soit une donnée de Shimura mixte.

On note $(\G_r,\X_r)$ la donnée de Shimura pure $(\QP_r,\Y_r)/\N_r$ (cf [P1] 2.9).
Elle est isomorphe à la donnée $(\GSp_{2r,\Q},\X_{r})$ de la section 1.1.

On a une application ``partie imaginaire'' de $\Y_r$ dans $\U_r(\R)(-1)=(2i\pi)^{-1}\U_r(\R)$ ([P1] 4.14),
et $\X^+$ (resp. $\X^-$) est l'image réciproque d'un cône ouvert convexe
$C(\X^+,\QP_r)$ (resp. $C(\X^-,\QP_r)$) de $\U_r(\R)(-1)$ ([P1] 4.15).
Identifions $\U_r(\R)(-1)$ à $M_{d-r}(\R)$.
Alors, d'après [P1] 4.26, le cône $C(\X^+,\QP_r)$ (resp. $C(\X^-,\QP_r)$) est l'ensemble des matrices
de la forme $AJ_{d-r}$, avec $A$ symétrique définie positive (resp. négative).

Si $\Pa$ est un sous-groupe parabolique maximal quelconque de $\G$,
il est conjugué par $\G(\Q)$ à l'un des $\Pa_r$,
et on peut donc associer à $\Pa$ une donnée de Shimura mixte $(\QP,\Y)$,
et une donnée de Shimura pure $(\G_P,\X_P)=(\QP,\Y)/\N$, 
où $\N$ est le radical unipotent de $\Pa$.
Les données de Shimura mixtes $(\QP,\Y)$ sont appelées \emph{composantes rationnelles de bord}
de $(\G,\X)$ ([P1] 4.11).

On fixe $n\geq 3$, et on note $\K=\K_d(n)$.

La compactification de Baily-Borel partielle de $\X$ est
$$\X^*=\X\sqcup\coprod_\Pa \X_P,$$
où $\Pa$ parcourt l'ensemble des sous-groupes paraboliques maximaux de $\G$;
on munit $\X^*$ de la topologie de Satake (cf par exemple [P1] 6.2).
L'action de $\G(\Q)$ sur $\X$ se prolonge en une action continue sur $\X^*$,
et la compactification de Baily-Borel de $M^{\K}(\G,\X)(\C)$ est 
$$M^{\K}(\G,\X)^*(\C)=\G(\Q)\sous (\X^*\times\G(\Af)/\K).$$
Elle a une structure canonique de variété complexe projective normale ([BB] 10.11).

Pink a montré que $M^{\K}(\G,\X)^*(\C)$ 
a un modèle canonique $M^{\K}(\G,\X)^*$ sur $\Q$ ([P1] 12.3),
qui est un schéma projectif normal sur $\Q$.
Le modèle canonique $M^{\K}(\G,\X)$ de $M^{\K}(\G,\X)(\C)$
est un ouvert dense de $M^{\K}(\G,\X)^*$.
De plus, toujours d'après [P1] 12.3, la stratification du bord de $M^{\K}(\G,\X)^*$ 
de [P1] 6.3 et [P2] 3.7 est définie sur $\Q$.

Rappelons la définition de cette stratification (on suit [P2] 3.7).
Soient $r\in\{0,\dots,d-1\}$ et $g\in\G(\Af)$. On pose
$\Hr_{g,r}=g\K g^{-1}\cap\Pa_r(\Q)\QP_r(\Af)$,
$\Hr_{\ell,g,r}=g\K g^{-1}\cap\Le_{\ell,r}(\Q)\N_r(\Af)$,
$\K_{Q,g,r}=g\K g^{-1}\cap\QP_r(\Af)$,
$\K_{N,g,r}=g\K g^{-1}\cap\N_r(\Af)$,
$\Gamma_{\ell,g,r}=\Hr_{\ell,g,r}/\K_{N,g,r}$
et $\K_{g,r}=\K_{Q,g,r}/\K_{N,g,r}\subset\G_r(\Af)$.
Le groupe $\Hr_{g,r}$ agit sur $M^{\K_{g,r}}(\G_r,\X_r)$
et son sous-groupe distingué d'indice fini $\Hr_{\ell,g,r}\K_{Q,g,r}$ agit trivialement.
On note $M^{\K_{g,r}}(\G_r,\X_r)/\Hr_{g,r}$ le quotient de la variété quasi-projective
$M^{\K_{g,r}}(\G_r,\X_r)$ par le groupe fini $\Hr_{g,r}/\Hr_{\ell,g,r}\K_{Q,g,r}$.
On a un morphisme
$$i_{g,r}:M^{\K_{g,r}}(\G_r,\X_r)\fl M^{\K}(\G,\X)^*$$
qui est le composé du morphisme fini 
$M^{\K_{g,r}}(\G_r,\X_r)\fl M^{\K_{g,r}}(\G_r,\X_r)/\Hr_{g,r}$
et d'une immersion localement fermée.
Sur les points complexes, $i_{g,r}$ est donné par le diagramme suivant:
$$\xymatrix{M^{K_{g,r}}(\G_r,\X_r)(\C)\ar@{=}[r] & \G_r(\Q)\sous (\X_r\times\G_r(\Af)/\K_{g,r}) & [(x,q\N_r(\Af))] \\
 & \QP_r(\Af)\sous (\X_r\times\QP_r(\Af)/\K_{Q,g,r})\ar[u]^{\wr}\ar[d] & [(x,q)]\ar@{|->}[u]\ar@{|->}[d] \\
M^{\K_{g,r}}(\G_r,\X_r)(\C)/\Hr_{g,r}\ar@{=}[r] & \Pa_r(\Q)\sous (\X_r\times\Pa_r(\Q)\QP_r(\Af)/\Hr_{g,r})\ar[d] & [(x,q)]\ar@{|->}[d] \\
M^{\K}(\G,\X)^*(\C)\ar@{=}[r] & \G(\Q)\sous (\X^*\times\G(\Af)/\K) & [(x,qg)]}$$ 
$i_{g,r}$ s'étend en un morphisme fini $\overline{i}_{g,r}:M^{\K_{g,r}}(\G_r,\X_r)^*\fl M^{\K}(\G,\X)^*$.
Les images des morphismes $i_{g,r}$ forment une stratification du bord de $M^{\K}(\G,\X)^*$,
et on a $Im(i_{g,r})=Im(i_{g',r'})$ si et seulement si $r=r'$ et
$\Pa_r(\Q)\QP_r(\Af)g\K=\Pa_r(\Q)\QP_r(\Af)g'\K$.

Nous allons voir que $i_{g,r}$ est une immersion.

\begin{lemme} Le double quotient $\Pa_r(\Q)\QP_r(\Af)\sous\G(\Af)/\G(\widehat{\Z})$ est un singleton.
\end{lemme}

Il existe donc un système de représentants $(g_j)_{j\in J}$ dans $\G(\widehat{\Z})$
du double quotient\newline
 $\Pa_r(\Q)\QP_r(\Af)\sous\G(\Af)/\K$.
Comme $\K$ est distingué dans $\G(\widehat{\Z})$, les groupes
$\Hr_{g_j,r}$, $\Hr_{\ell,g_j,r}$, etc, ne dépendent pas de $j$;
on les note donc $\Hr_r$, $\Hr_{\ell,r}$, etc.
Comme $\K=Ker(\G(\widehat{\Z})\fl\G(\Z/n\Z))$,
on voit que $\Hr_r/\K_{N,r}=\Gamma_{\ell,r}\times\K_r$, avec
$\Gamma_{\ell,r}=Ker(\Le_{\ell,r}(\Z)\fl\Le_{\ell,r}(\Z/n\Z))$ et
$\K_r=Ker(\G_r(\widehat{\Z})\fl \G_r(\Z/n\Z))$.
En particulier, $\Hr_r=\Hr_{\ell,r}\K_{Q,r}$, 
et les $i_{g_j,r}$, $j\in J$, sont des immersions localement fermées.
De plus, les strates $M^{\K_r}(\G_r,\X_r)$ sont isomorphes à la variété de Shimura
$\Mod_{r,n,\Q}$ de la section 1.1.

La compactification de Baily-Borel $M^{\K}(\G,\X)^*$ n'est pas lisse si $d\geq 3$,
mais elle admet une famille de résolutions des singularités, 
les compactifications toroïdales (cf [AMRT] ou [P1] chapitres 6-9 pour la construction sur $\C$,
[P1] chapitre 12 pour les modèles canoniques).
Ces compactifications dépendent d'une décomposition en cônes admissible pour $(\G,\X)$ ([P1] 6.4).
Nous allons définir certaines décompositions particulières.

Il faut d'abord rappeler quelques définitions. 

Soit $(\QP,\Y)$ une composante rationnelle de bord de $(\G,\X)$. 
On note (cf [P1] 4.22) $C^*(\X^+,\QP)$ (resp. $C^*(\X^-,\QP)$) l'union des cônes
$C(\X^+,\QP')$ (resp. $C(\X^-,\QP')$) pour $(\QP',\Y')$ parcourant
l'ensemble des composantes rationnelles de bord entre $(\QP,\Y)$ et $(\G,\X)$
(ie telles que $\QP\subset\QP'$).
D'après [P1] 4.26, pour tout $r\in\{0,\dots,d-1\}$,
$C^*(\X^+,\QP_r)\subset\U_r(\R)(-1)$ (resp. $C^*(\X^-,\QP_r)\subset\U_r(\R)(-1)$)
est égal à l'ensemble des matrices $AJ_{d-r}$, où $A$ est symétrique positive (resp. négative)
et le noyau de $A$ est défini sur $\Q$.

Le complexe conique $\mathcal{C}(\G,\X)$ de $(\G,\X)$ est
le quotient de l'union disjointe sur l'ensemble des composantes rationnelles de bord $(\QP,\Y)$
des $C^*(\X^+,\QP)\sqcup C^*(\X^-,\QP)$
par la relation d'équivalence engendrée par les graphes des inclusions
$C^*(\X^+,\QP')\sqcup C^*(\X^-,\QP')\subset C^*(\X^+,\QP)\sqcup C^*(\X^-,\QP)$
pour $\QP\subset\QP'$
([P1] 4.24).

Une décomposition en cônes $\K$-admissible $\Sc_{ad}$ pour $(\G,\X)$
est une collection de sous-ensembles de $\mathcal{C}(\G,\X)\times\G(\Af)$
qui vérifie les conditions (i) à (v) de [P1] 6.4.
En particulier, on demande que pour toute composante rationnelle de bord
$(\QP,\Y)$, pour toute composante connexe $\X^\circ$ de $\X$ et pour tout $g\in\G(\Af)$,
$$\Sc_{ad}(\X^\circ,\QP,g)=\{\sigma\in\Sc_{ad},\sigma\subset C^*(\X^\circ,\QP)\times\{g\}\}$$
soit une décomposition (partielle) en cônes polyédraux rationnels du cône 
$C^*(\X^\circ,\QP)\times\{g\}$.

On note $C$ le cône convexe de $M_d(\R)$ des matrices symétriques positives
dont le noyau est défini sur $\Q$; on fait agir $\GL_d(\Z)$ sur $C$ par 
$(g,A)\fle gA{}^tg$.
Soit $\Sc$ une décomposition en cônes polyédraux de $C$ 
invariante par l'action de $\GL_d(\Z)$ et telle que
$\GL_d(\Z)\sous\Sc$ soit fini.
On veut lui associer une décomposition en cônes admissible $\Sc_{ad}$ pour $(\G,\X)$.

Il suffit de définir $\Sc_{ad}(\X^{\pm},\QP,g)$ pour les composantes de bord minimales $(\QP,\Y)$ de $(\G,\X)$.
Soit $(\QP,\Y)$ une telle composante. On se ramène par conjugaison au cas où $(\QP,\Y)=(\QP_0,\Y_0)$.
On a vu plus haut que le cône $C^*(\X^+,\QP_0)$ (resp. $C^*(\X^-,\QP_0)$) était égal à $CJ_d$ 
(resp. $-CJ_d$).
Pour tout $g\in\G(\Af)$, on pose
$$\Sc_{ad}(\X^+,\QP_0,g)=\{\sigma.J_d,\sigma\in\Sc\}\times\{g\}$$
$$\Sc_{ad}(\X^-,\QP_0,g)=\{-\sigma.J_d,\sigma\in\Sc\}\times\{g\}.$$
Il est clair qu'on obtient bien une décomposition en cônes admissible pour $(\G,\X)$.
De plus, cette décomposition est complète (resp. lisse pour $\K$)
si et seulement si $\Sc$ est complète (resp. lisse pour le réseau $nM_d(\Z)\subset M_d(\R)$, cf [P1] 5.2).

On note $M^{\K}(\G,\X,\Sc)$ la compactification toroïdale associée.
Si $\Sc$ est complète et lisse, $M^{\K}(\G,\X,\Sc)$ est projective et le bord
$M^{\K}(\G,\X,\Sc)-M^{\K}(\G,\X)$ est une union de diviseurs lisses à croisements normaux.

L'identité $M^{\K}(\G,\X)$ s'étend en un morphisme surjectif $M^{\K}(\G,\X,\Sc)\fl M^{\K}(\G,\X)^*$.
Les images inverses des strates de $M^{\K}(\G,\X)^*$ forment une stratification du bord de $M^{\K}(\G,\X,\Sc)$,
et on sait décrire ces strates et les complétés formels de $M^{\K}(\G,\X,\Sc)$ le long de ces strates
en termes de certains plongements toriques ([P2] 3.10).

\subsection{Compactifications sur $\Z$}
\hspace{.5cm}

Dans cette section, $(\G,\X)$ est la donnée de Shimura $(\GSp_{2d,\Q},\X_d)$ de la section 1.1.
On suppose toujours $n\geq 3$, et on note $\K=\K(n)$.

Soit $\Sc$ une décomposition en cônes polyédraux du cône $C$ de la section 1.2.
On suppose que $\Sc$ est invariante par l'action de $\GL_d(\Z)$, complète et lisse
pour le réseau $M_d(\Z)$,
et que $\GL_d(\Z)\sous\Sc$ est fini.
Chai et Faltings ([CF] IV.6.7) ont montré que 
la compactification toroïdale $M^{\K}(\G,\X,\Sc)$ a un modèle entier $\Mod_{d,n}(\Sc)\supset\Mod_{d,n}$ :
$\Mod_{d,n}(\Sc)$ est un espace algébrique propre et lisse sur $\Z[1/n]$
dont la fibre générique est $M^{\K}(\G,\X,\Sc)$,
et le bord $\Mod_{d,n}(\Sc)-\Mod_{d,n}$ est un diviseur à croisements normaux relatif sur $\Z[1/n]$.
Si $n$ divise $m$, le morphisme $T_1:\Mod_{d,m}\fl\Mod_{d,n}$ se prolonge en un morphisme
$\Mod_{d,m}(\Sc)\fl\Mod_{d,n}(\Sc)$, qu'on notera $\til{T}_1$.

En utilisant les compactifications toroïdales, Chai et Faltings ont aussi construit 
un modèle entier de $M^{\K}(\G,\X)^*$ ([CF] V.2.5) :
un schéma $\Mod_{d,n}^*\supset\Mod_{d,n}$ normal et projectif sur $\Z[1/n]$,
dont la fibre générique est $M^{\K}(\G,\X)^*$.
De plus, on a une stratification de $\Mod_{d,n}^*-\Mod_{d,n}$ par des
$\Mod_{r,n}$, $0\leq r\leq d-1$, 
qui étend la stratification de la section 1.2,
et l'adhérence d'une strate est isomorphe à sa compactification de Baily-Borel ([CF] V.2.5 (4)).
On notera toujours $i_{g,r}$ et $\overline{i}_{g,r}$ les prolongements des morphismes
$i_{g,r}$ et $\overline{i}_{g,r}$ de 1.2.
L'action de $\GSp_{2d}(\Z/n\Z)$ sur $\Mod_{d,n}$ s'étend à $\Mod^*_{d,n}$ ([CF] V.2.5 (3)),
donc les morphismes $T_g$ définis dans la section 1.1 se prolongent
aux compactifications de Baily-Borel.
Nous noterons $\overline{T}_g$ ces prolongements.

\begin{lemme} Soient $n,m\geq 3$ deux entiers tels que $n$ divise $m$ et $p$ un nombre premier
qui ne divise pas $m$. Alors, pour tout $g\in\G(\Af^p)$ tel que
$g^{-1}\K_d(m)g\subset\K_d(n)$, le morphisme
$T_g:\Mod_{d,m|\Z_{(p)}}\fl\Mod_{d,n|\Z_{(p)}}$ de 1.1 se prolonge en un morphisme
$\overline{T}_g:\Mod^*_{d,m|\Z_{(p)}}\fl\Mod^*_{d,n|\Z_{(p)}}$.

\end{lemme}

\begin{proof} Supposons d'abord $d\geq 2$. On note $\omega$ le faisceau canonique sur
$\Mod_{d,n}$ et sur $\Mod_{d,m}$ (cf \cite{CF} I.4.11). 
D'après le principe de Koecher (\cite{CF} V.1.8 (iii)) et
le point (3) du théorème V.2.5 de \cite{CF}, le schéma $\Mod_{d,n}^*$ est canoniquement
isomorphe à $\Proj(\bigoplus\limits_{k\in\Nat}\Gamma(\Mod_{d,n},\omega^{\otimes k}))$,
et on a une formule analogue pour $\Mod_{d,m}^*$. Le lemme en résulte.

Si $d=1$, il n'existe qu'une seule décomposition en cônes polyédraux complète $\Sc$, et
on a $\Mod_{d,n}^*=\Mod_{d,n}(\Sc)$, $\Mod_{d,m}^*=\Mod_{d,m}(\Sc)$. Le lemme est une
conséquence du corollaire V.6.11 de \cite{CF}.

Si $d=0$, on a $\Mod_{d,n}=\Mod_{d,n}^*$ et $\Mod_{d,m}=\Mod_{d,m}^*$, donc il n'y a rien
à prouver.

\end{proof}

Enfin, le morphisme $M^{\K}(\G,\X,\Sc)\fl M^{\K}(\G,\X)^*$ se prolonge en un morphisme
$\Mod_{d,n}(\Sc)\fl\Mod_{d,n}^*$ ([CF] IV.2.5 (2)),
et on a une description analogue à celle de [P2] 3.10 des images inverses des strates de $\Mod_{d,n}^*$ ([CF] IV.6.7 (4), IV.6.11 et V.2.5 (5)).

\section{Systèmes de coefficients et théorème de Pink}

Dans cette section, $(\G,\X)$ est la donnée de Shimura $(\GSp_{2d,\Q},\X_d)$ 
de la section 1.1,
$n$ est un entier $\geq 3$, $\K=\K_d(n)\subset\G(\Af)$,
et $\ell$ est un nombre premier.

\subsection{Systèmes de coefficients}
\hspace{.5cm}

Si $V$ est une représentation algébrique de $\G(\Q_\ell)$ dans un $\Q_\ell$-espace vectoriel de dimension finie,
on peut lui associer un système local de $\Q_\ell$-espaces vectoriels sur $M^{\K}(\G,\X)(\C)$,
$$\G(\Q)\sous (V\times\X\times\G(\Af)/\K)\fl \G(\Q)\sous(\X\times\G(\Af)/\K)=M^{\K}(\G,\X)(\C)$$
($\G(\Q)$ agit diagonalement).
Langlands a montré que ces systèmes de coefficients provenaient de faisceaux $\ell$-adiques lisses
sur $M^{\K}(\G,\X)_\C$ (cf [L] p 34-38).
Nous allons rappeler ici la méthode de Pink pour étendre ces systèmes de coefficients au modèle entier.

Soit $\varphi:\widehat{X}\fl X$ un revêtement étale galoisien de groupe profini $\Gamma$.
Pink a construit dans [P1] 1.10 un foncteur exact
$$\mu_{\Gamma,\varphi}:Mod_\Gamma\fl {\acute E}t_X$$
de la catégorie $Mod_\Gamma$ des $\Gamma$-modules à gauche continus (discrets)
dans la catégorie ${\acute E}t_X$ des faisceaux abéliens étales sur $X$.
Pour tout $M\in Mod_\Gamma$, on a
$$\mu_{\Gamma,\varphi}M=(M\otimes\limid \varphi_{\Delta*}\underline{\Z})^\Gamma,$$
où $\Delta$ parcourt l'ensemble des sous-groupes ouverts distingués de $\Gamma$ et
$\varphi_\Delta$ est le revêtement étale fini
$\widehat{X}/\Delta\fl X$.
Les fibres de $\mu_{\Gamma,\varphi}M$ sont données par la proposition suivante :

\begin{proposition}\label{fibres_faisceaux}([P2] 1.10.4) Soit $M$ un objet de $Mod_\Gamma$. 
Soit $x_0$ un point de $X$,
$k$ son corps résiduel, 
$\overline{k}$ une clôture séparable de $k$, 
$x:Spec(\overline{k})\fl X$ le point géométrique de $X$ correspondant, 
$\widehat{x}:Spec(\overline{k})\fl\widehat{X}$ un point géométrique de $\widehat{X}$ au-dessus de $x$. 
Pour tout $\sigma\in Gal(\overline{k}/k)$, on note $\psi(\sigma)$ l'unique élément de $\Gamma$
 tel que $\sigma.\widehat{x}=\widehat{x}.\psi(\sigma)$. Alors 
$\psi:Gal(\overline{k}/k)\fl\Gamma$
est un morphisme continu, 
et la fibre de $\mu_{\Gamma,\varphi}(M)$ en $x$ est isomorphe à $M$ 
avec l'action de $Gal(\overline{k}/k)$ donnée par $\sigma.m=\psi(\sigma).m$. 

\end{proposition}

De plus, on sait calculer les images inverses des $\mu_{\Gamma,\varphi}M$
par certains revêtements étales, grâce au lemme ci-dessous, qui est un cas particulier de [P2] 1.11.5 :

\begin{lemme}\label{lemme:Shapiro} Soit $\Gamma'$ un sous-groupe fermé de $\Gamma$.
On note $\varphi':\widehat{X}\fl X'=\widehat{X}/\Gamma'$ 
et $f:X'\fl X$ les morphismes évidents.

Alors, pour tout $M\in Mod_\Gamma$, on a un isomorphisme canonique
$$f^*\mu_{\Gamma,\varphi}M\simeq\mu_{\Gamma',\varphi'}Res_{\Gamma'}^{\Gamma}M.$$

\end{lemme}

Revenons à la situation du début.
On note $Rep_{\G(\Q_\ell)}$ la catégorie des représentations algébriques de $\G(\Q_\ell)$
dans des $\Q_\ell$-espaces vectoriels de dimension finie,
et $\K_\ell$ l'image de $\K$ dans $\G(\Q_\ell)$.
On fait agir $\G(\Af)$ sur les objets de $\G(\Q_\ell)$ via
la projection $\G(\Af)\fl\G(\Q_\ell)$.

Soit
$$\widehat{\Mod}_{d,n}=\lim_{\overleftarrow{r\in\Nat}}\Mod_{d,\ell^rn}.$$
$\widehat{\Mod}_{d,n}$ est un schéma sur $\Z[1/\ell n]$ (il n'est pas localement de type fini),
et le morphisme évident $\varphi:\widehat{\Mod}_{d,n}\fl\Mod_{d,n}[1/\ell]=\Mod_{d,n|Spec(\Z[1/n\ell])}$
est un revêtement étale galoisien de groupe $\K_\ell$.

\begin{definition}\label{def_faisceaux} Soit $V\in Rep_{\G(\Q_\ell)}$.
On choisit un $\Z_\ell$-réseau $\Lambda\subset V$ invariant par $\K_\ell$,
et on pose
$$\F^{\K}V=\Q_\ell\otimes\lim_{\overleftarrow{m\in\Nat}}\mu_{\K_\ell,\varphi}(\Lambda/\ell^m\Lambda).$$
C'est un faisceau $\ell$-adique lisse sur $\Mod_{d,n}[1/\ell]$,
qui à isomorphisme près ne dépend pas du choix de $\Lambda$.
On obtient donc un foncteur exact $\F^{\K}$ de $Rep_{\G(\Q_\ell)}$ dans la
catégorie des faisceaux $\ell$-adiques sur $\Mod_{d,n}[1/\ell]$,
et on notera encore $\F^{\K}$ son foncteur dérivé.

\end{definition}

Pour toute $V\in Rep_{\G(\Q_\ell)}$, l'analytisé de 
l'image réciproque de $\F^{\K}V$ sur $M^{\K}(\G,\X)_\C$
est le système de coefficients défini au début de cette section (cf [P2] 5.1).

On s'intéresse enfin aux poids des faisceaux obtenus. 
On dira qu'une représentation $V$ de $\G(\Q_\ell)$ est pure de poids $t$
si le centre $\Gr_{m,\Q}.I_{2d}$ de $\G$ 
agit sur $V$ par le caractère $x\fle x^t$,
et qu'elle est de poids $<t$ (resp. $\geq t$) si elle est somme de sous-représentations
pures dont les poids $t'$ vérifient $t'<t$ (resp. $t'\geq t$).

\begin{proposition}\label{poids_faisceaux}([P2] 5.6.6, voir aussi [LR2] 6) Pour toute représentation $V\in Rep_{\G(\Q_\ell)}$, le faisceau $\ell$-adique
lisse $\F^{\K}V$ est mixte (au sens de [D] 1.2.2),
et il est pur de poids $-t$ si $V$ est pure de poids $t$.

\end{proposition}

\subsection{Prolongements des systèmes de coefficients à la compactification de Baily-Borel}
\hspace{.5cm}

Dans cette section, on utilise les notations des sections 1.2 et 2.1.
En particulier, pour tout $r\in\{0,\dots,d-1\}$, on note 
$\Hr_{\ell,r}=\K\cap\Pa_{\ell,r}(\Q)\N_r(\Af)$ et
$\Gamma_{\ell,r}=\Hr_{\ell,r}/(\Hr_{\ell,r}\cap\N_r(\Af))$
(c'est un sous-groupe arithmétique de la partie linéaire $\Le_{\ell,r}$ de $\Le_r$).

Nous allons énoncer le théorème principal de [P2] pour les variétés de Shimura considérées.
Remarquons d'abord que, comme expliqué dans [P2] 4.9, comme on dispose des modèles entiers $\Mod_{d,\ell^rn}$, $r\in\Nat$,
et de compactifications de Baily-Borel et toroïdales de ces modèles sur $\Z[1/\ell n]$ 
qui vérifient les propriétés de [P2] 3.7-3.11,
la preuve du théorème [P2] 4.2.1 s'étend aux modèles entiers si on ne considère que des
faisceaux de $\Z_\ell$-torsion.
On obtient donc le théorème suivant :

\begin{theoreme}\label{th_Pink} On note $j$ l'immersion ouverte $\Mod_{d,n}[1/\ell]\fl\Mod_{d,n}^*[1/\ell]$.
Pour tout $V\in D^b(Rep_{\G(\Q_\ell)})$, pour tous $r\in\{0,\dots,d-1\}$ et
$g\in\G(\widehat{\Z})$, on a un isomorphisme canonique
$$\begin{array}{rcl}i_{g,r}^*Rj_*\F^{\K}V & \simeq & \F^{\K_r}R\Gamma(\Hr_{\ell,r},g.V) \\
& \simeq & \F^{\K_r}R\Gamma(\Gamma_{\ell,r},R\Gamma(Lie(\N_r),g.V)),\end{array}$$
où $g.V$ est $V$ avec l'action de $\G(\Q_\ell)$ donnée par $(h,v)\fle (g_\ell^{-1}hg_\ell).v$.

\end{theoreme}

\begin{proof} Le premier isomorphisme est celui de [P2] 4.2.1,
compte tenu des remarques ci-dessus et 
du fait que $\Hr_r/\Hr_{\ell,r}=\K_r$.
Le deuxième isomorphisme vient de l'isomorphisme canonique $R\Gamma(\Hr_{\ell,r},\ )\simeq R\Gamma(\Gamma_{\ell,r},R\Gamma(\K_{N,r},\ ))$,
du lemme 5.2.2 de [P2]
(qui dit que la cohomologie continue de $\K_{N,r}$ est égale à la cohomologie du
groupe discret $\K_{N,r}\cap\N_r(\Q)$)
et du lemme ci-dessous, 
qui est une version algébrique du théorème de van Est.

\end{proof}

\begin{lemme}\label{van:Est}Soient $\U$ un groupe algébrique unipotent connexe sur $\Q$ et $\Gamma_U$ un sous-groupe arithmétique de $\U(\Q)$. 
On note $Mod_\U$ la catégorie des limites inductives de représentations algébriques rationnelles de dimension finie de $\U$ (ou, ce qui revient au même, de $Lie(\U)$),
et $A_U$ la $\Q$-algèbre des polynômes à coefficients rationnels sur $\U(\Q)$, qui contient la sous-algèbre $\Q$ des constantes.
Pour tout $M$, on définit une suite exacte
$$0\fl M\stackrel{u_0}{\fl}I^0(M)\stackrel{u_1}{\fl} I^1(M)\stackrel{u_2}{\fl}I^2(M)\stackrel{u_3}{\fl}\dots$$
par :
\begin{itemize}
\item[$\bullet$] $u_0:M\fl I^0(M)$ est l'injection évidente $M\simeq M\otimes\Q\fl M\otimes A_U=I^0(M)$;
\item[$\bullet$] pour tout $i\in\Nat$, $u_{i+1}:I^i(M)\fl I^{i+1}(M)$ est le morphisme évident $I^i(M)\fl Coker(u_i)\fl Coker(u_i)\otimes A_U=I^{i+1}(M)$.

\end{itemize}

Alors :
\begin{itemize}
\item[(i)] Pour tout $M\in Mod_\U$, $M\otimes A_U$ est un objet injectif de $Mod_\U$ et un objet acyclique pour le foncteur $(\quad)^{\Gamma_U}$ (invariants par $\Gamma_U$).
\item[(ii)] Pour tout $M\in Mod_\U$, le morphisme suivant est un isomorphisme :
$$R\Gamma(Lie(\U),M)=R\Gamma(\U(\Q),M)\simeq I^\bullet(M)^{\U(\Q)}\fl I^{\bullet}(M)^{\Gamma_U}\simeq R\Gamma(\Gamma_U,M).$$

\end{itemize}

\end{lemme}

Ce lemme est prouvé dans [GHM] 24.

Enfin, la compatibilité de l'isomorphisme du théorème de Pink avec
les morphismes $T_g$ de la section 1.1 est explicitée dans la proposition ci-dessous :

\begin{proposition}\label{Pink:corr_Hecke}(\cite{P2} 4.8.5) On suppose que $n$ divise $m$,
on fixe un nombre premier $p\not=\ell$ qui ne divise pas $m$ 
et on note $\K'=\K_d(m)$ et 
$j':\Mod_{d,m}[1/\ell]\fl\Mod^*_{d,m}[1/\ell]$ l'inclusion.
Soient $g\in\G(\Af^p)$ tel que
$g^{-1}\K' g\subset\K$, $h\in\G(\widehat{\Z})$ et $r\in\{0,\dots,d-1\}$. On fixe
$h'\in\G(\widehat{\Z})$, $q\in\P_r(\Q)\QP_r(\Af)\cap\G(\Af^p)$ et $k\in\K$ tels que
$hg=qh'k$.
On note $i'_{h,r}$ l'inclusion de la strate $\Mod_{r,m}$ dans $\Mod^*_{d,m}$ définie par $h$.

Alors le morphisme $\overline{T}_g:\Mod^*_{d,m|\Z_{(p)}}\fl\Mod^*_{d,n|\Z_{(p)}}$
envoie l'image de $i'_{h,r}$ dans celle de $i_{h',r}$, et
le morphisme $\Mod_{r,m|\Z_{(p)}}\fl\Mod_{r,n|\Z_{(p)}}$ obtenu en restreignant
$\overline{T}_g$ à ces images est $T_{\overline{q}}$, où $\overline{q}$ est l'image de
$q$ dans $\G_r(\Af^p)$.

Soit $V\in D^b(Rep_{\G(\Q_\ell)})$.
D'après le théorème de Pink et le lemme \ref{lemme:Shapiro}, 
on a des isomorphismes canoniques
$$T_{\overline{q}}^*i_{h',r}^*Rj_*\F^{\K}V\simeq T_{\overline{q}}^*\F^{\K_r}R\Gamma
(\Hr_{\ell,r},h'.V)\simeq \F^{\K'_r}R\Gamma(\Hr_{\ell,r},qh'.V)\simeq
\F^{\K'_r}R\Gamma(\Hr_{\ell,r},hg.V)$$
$${i'_{h,r}}^*Rj'_*T_g^*\F^{\K}V\simeq {i'_{h,r}}^*Rj'_*\F^{\K'}g.V\simeq \F^{\K'_r}R\Gamma(\Hr'_{\ell,r},hg.V),$$
où le dernier isomorphisme de la première ligne vient de l'isomorphisme
$qh'.V\iso qh'k.V=hg.V$, $v\fle k.v$.
Par ces isomorphismes, le morphisme de changement de base
$$T_{\overline{q}}^*i_{h',r}^*Rj_*\F^{\K}V={i'_{h,r}}^*\overline{T}_g^*Rj_*\F^{\K}(V)
\stackrel{CB}{\fl}{i'_{h,r}}^*Rj'_*T_g^*\F^{\K}(V)$$
correspond au morphisme induit par l'inclusion $\Hr'_{\lin,r}\subset\Hr_{\lin,r}$
$$\F^{\K'_r}R\Gamma(\Hr_{\lin,r},hg.V)\fl\F^{\K'_r}R\Gamma(\Hr'_{\lin,r},hg.V).$$

\end{proposition}

\section{Prolongement intermédiaire des faisceaux pervers purs}

Cette partie est indépendante des autres.

On fixe un corps fini $\Fi_q$ et un nombre premier $\ell$ inversible dans $\Fi_q$.
Tous les schémas sont séparés de type fini sur $\Fi_q$.
Si $X$ est un schéma, on note $D^b_m(X,\Q_\ell)$ la catégorie des complexes $\ell$-adiques mixtes sur $X$ (au sens de [BBD] 5.1.5; en particulier, les complexes sont à poids entiers),
munie de la t-structure donnée par la perversité autoduale.

Soient $X$ un schéma et $j:U\fl X$ un ouvert non vide de $X$.
Le premier objectif de cette partie est d'écrire $j_{!*}K$, où $K$ est un faisceau pervers pur sur $U$,
comme un certain tronqué par le poids de $Rj_*K$ (théorème \ref{th:calcul_IC_general}).
Le deuxième objectif est
 de calculer la trace d'une puissance $\Phi$ de l'endomorphisme de Frobenius
sur la cohomologie de $j_{!*}K$ en fonction de
la trace de $\Phi$ sur la cohomologie de $Rj_*K$ et de complexes 
de même type supportés par $X-U$
(théorème \ref{th:simplification_formule_traces}).

\subsection{Troncature par le poids dans $D^b_m(X,\Q_\ell)$}
\hspace{.5cm}

\begin{proposition}\label{t_structures:prop1} Soit $X$ un schéma sur $\Fi_q$.
Pour tout $a\in\Z\cup\{\pm\infty\}$, on note $\DP^{\leq a}(X)$, ou $\DP^{\leq a}$ s'il n'y a pas d'ambiguïté sur $X$, (resp. $\DP^{\geq a}(X)$ ou $\DP^{\geq a}$)
la sous-catégorie pleine de $D^b_m(X,\Q_\ell)$ dont les objets sont les complexes mixtes $K$ 
tels que pour tout $i\in\Z$, $\Hp^iK$ soit de poids $\leq a$ (resp. $\geq a$).
Alors:
\begin{itemize}
\item[(i)] $\DP^{\leq a}$ et $\DP^{\geq a}$ sont des sous-catégories stables par décalage et extensions (cf [BBD] 1.2.6) de $D^b_m(X,\Q_\ell)$ (en particulier, ce sont des sous-catégories triangulées).\newline
La dualité de Poincaré échange $\DP^{\leq a}$ et $\DP^{\geq -a}$.\newline
Si $a<a'$, on a $\DP^{\leq a}\cap \DP^{\geq a'}=0$.
\item[(ii)] On a $\DP^{\leq a}(1)=\DP^{\leq a-2}$ et $\DP^{\geq a}(1)=\DP^{\geq a-2}$
(où $(1)$ est le twist à la Tate).
\item[(iii)] Pour tous $K\in\DP^{\leq a}$ et $L\in\DP^{\geq a+1}$, on a
$R\Hom(K,L)=0$.
\item[(iv)] Pour tout $a\in\Z\cup\{\pm\infty\}$, $(\DP^{\leq a},\DP^{\geq a+1})$ est une t-structure sur $D^b_m(X,\Q_\ell)$.
\end{itemize}
\end{proposition}

Cette proposition sera démontrée dans la section 3.2.

D'après [BBD] 1.3.3,
l'inclusion $\DP^{\leq a}\subset D^b_m(X,\Q_\ell)$ (resp. $\DP^{\geq a}\subset D^b_m(X,\Q_\ell)$) 
admet un adjoint à droite (resp. à gauche), qu'on notera $w_{\leq a}$ (resp. $w_{\geq a}$),
et pour tout $K\in D^b_m(X,\Q_\ell)$, il existe un unique morphisme $w_{\geq a+1}K\fl (w_{\leq a})K[1]$ qui fait du triangle suivant un triangle distingué
$$w_{\leq a}K\fl K\fl w_{\geq a+1}K\fl (w_{\leq a}K)[1].$$
Ce triangle est, à isomorphisme unique près, l'unique triangle distingué $A\fl K\fl B\stackrel{+1}{\fl}$ avec $A\in\DP^{\leq a}$ et $B\in\DP^{\geq a+1}$ (toujours d'après [BBD] 1.3.3).

Comme la dualité de Poincaré échange $\DP^{\leq a}$ et $\DP^{\geq -a}$, 
elle échange aussi $w_{\leq a}$ et $w_{\geq -a}$.

\begin{remarques}
\begin{itemize}
\item[(1)] $\DP^{\leq a}$ n'est pas la catégorie des complexes mixtes de poids $\leq a$ :
un complexe mixte $K$ est de poids $\leq a$ si et seulement si $\Hp^iK$ est de poids $\leq a+i$ pour tout $i\in\Z$ ([BBD] 5.4.1),
et cette condition est différente de la condition qui caractérise les objets de $\DP^{\leq a}$.
\item[(2)] Le décalage de $1$ dans la définition de la t-structure est nécessaire : 
$(\DP^{\leq a},\DP^{\geq a})$ n'est pas une t-structure.
\item[(3)] $(\DP^{\leq a},\DP^{\geq a+1})$ est une t-structure quelque peu étrange : 
son coeur est nul, et elle ne donne donc pas lieu à une théorie cohomologique intéressante.
De plus, $\DP^{\leq a}$ et $\DP^{\geq a+1}$ sont des sous-catégories triangulées de $D^b_m(X,\Q_\ell)$
(en particulier, elles sont stables par le foncteur $[1]$),
ce qui est inhabituel.
\item[(4)] Si $K$ est un faisceau pervers mixte, $w_{\leq a}$ est simplement le plus grand sous-faisceau pervers de $K$ de poids $\leq a$,
et $w_{\geq a}K$ est le plus grand quotient pervers de $K$ de poids $\geq a$ ([BBD] 5.3.5).\newline
Beilinson a montré dans [B] que le foncteur ``réalisation'' ([BBD] 3.1.9) de la catégorie dérivée
bornée de la catégorie des faisceaux pervers mixtes sur $X$ dans $D^b_m(X,\Q_\ell)$
est une équivalence de catégories.
Si $K\in D^ b_m(X,\Q_\ell)$ est représenté par un complexe $C^\bullet$ de faisceaux 
pervers mixtes, alors $w_{\leq a}K$ (resp. $w_{\geq a}K$) est représenté par le complexe
$(w_{\leq a}C^n)$ (resp. $(w_{\geq a}C^n)$).

\end{itemize}
\end{remarques}

La proposition suivante donne quelques propriétés des foncteurs $w_{\leq a}$ et $w_{\geq a}$.

\begin{proposition}\label{t_structures:prop2}
\begin{itemize}
\item[(i)] Pour tout $K\in D^b_m(X,\Q_\ell)$, on a $w_{\leq a}(K(1))=(w_{\leq a+2}K)(1)$ et $w_{> a}(K(1))=(w_{> a+2}K)(1)$.
\item[(ii)] Soit $K\in D^b_m(X,\Q_\ell)$. 
L'image par $\Hp^i$ de la flèche de cobord $w_{\geq a+1}K\fl (w_{\leq a}K)[1]$ est nulle pour tout $i\in\Z$,
donc la suite exacte longue de cohomologie perverse du triangle distingué
$$w_{\leq a}K\fl K\fl w_{\geq a+1}K\stackrel{+1}{\fl}$$
donne des suites exactes courtes de faisceaux pervers
$$0\fl\Hp^iw_{\leq a}K\fl\Hp^iK\fl \Hp^iw_{\geq a+1}K\fl 0.$$
\item[(iii)] $w_{\leq a}$ et $w_{\geq a}$ commutent au foncteur de décalage $[1]$,
et ils envoient les triangles distingués
de $D^b_m(X,\Q_\ell)$ sur des triangles distingués.
\item[(iv)] $w_{\leq a}$ et $w_{\geq a}$ envoient la catégorie abélienne des faisceaux pervers mixtes dans elle-même, 
et leurs restrictions à cette catégorie sont des foncteurs exacts.
Pour tout $K\in D^b_m(X,\Q_\ell)$, pour tout $i\in\Z$, on a
$$w_{\leq a}(\Hp^iK)=\Hp^iw_{\leq a}K$$
$$w_{\geq a}(\Hp^iK)=\Hp^iw_{\geq a}K.$$
\item[(v)] Soit $f:X\fl Y$ un morphisme. 
Si la dimension des fibres de $f$ est inférieure ou égale à $d$, alors
$$Rf_!(\DP^{\leq a}(X))\subset \DP^{\leq a+d}(Y)$$
$$Rf_*(\DP^{\geq a}(X))\subset \DP^{\geq a-d}(Y)$$
$$f^*(\DP^{\leq a}(Y))\subset \DP^{\leq a+d}(X)$$
$$f^!(\DP^{\geq a}(Y))\subset\DP^{\geq a-d}(X).$$

\end{itemize}
\end{proposition}

La démonstration des propositions \ref{t_structures:prop1} et \ref{t_structures:prop2} sera donnée dans la section 3.2.

\begin{theoreme}\label{th:calcul_IC_general} Soient $a\in\Z$, $X$ un schéma séparé 
de type fini sur $\Fi_q$, $j:U\fl X$ un ouvert non vide de $X$,
et $K$ un faisceau pervers pur de poids $a$ sur $U$.
Alors les flèches canoniques
$$w_{\geq a}j_!K\fl j_{!*}K\fl w_{\leq a}Rj_*K$$
sont des isomorphismes.

\end{theoreme}

Si $K$ est le faisceau constant $\Q_\ell$, cette formule est à rapprocher des formules 4.5.7 et 4.5.9
de l'article [S] de Morihiko Saito.

\begin{proof} Il suffit de montrer que la deuxième flèche est un isomorphisme 
(le cas de la première flèche en résulte par dualité).\newline
Cela découle des trois points suivants ($i$ est l'inclusion de $X-U$ dans $X$) :
\begin{itemize}
\item[(1)] Un complexe $K\in D^b_m(U,\Q_\ell)$ a au plus un prolongement $L\in D^b_m(X,\Q_\ell)$ tel que
$i^*L\in\DP^{\leq a}$ et $i^!L\in\DP^{\geq a+1}$.\newline
En effet, soient $L,L'\in D^b_m(X,\Q_\ell)$.
On a un triangle distingué
$$R\Hom(i^*L,i^!L')\fl R\Hom(L,L')\fl R\Hom(j^*L,j^*L')\stackrel{+1}{\fl}.$$
Si $i^*L\in\DP^{\leq a}$ et $i^!L\in\DP^{\geq a+1}$, alors $R\Hom(i^*L,i^!L')=0$ d'après le (iii) de la proposition \ref{t_structures:prop1},
donc on a un isomorphisme
$$R\Hom(L,L')\iso R\Hom(j^*L,j^*L').$$
\item[(2)] Soit $K\in D^b_m(U,\Q_\ell)$. On note $L=w_{\leq a}Rj_*K$. 
Alors $i^*L\in\DP^{\leq a}$ par le (iv) de la proposition \ref{t_structures:prop2}.
De plus, on a un triangle distingué
$$L\fl Rj_*K\fl w_{\geq a+1}Rj_*K\stackrel{+1}{\fl},$$
d'où un isomorphisme
$$i^!w_{\geq a+1}Rj_*K[-1]\iso i^!L.$$
D'après le (iv) de la proposition \ref{t_structures:prop2}, $i^!L\in\DP^{\geq a+1}$.
\item[(3)] Soit $K$ un faisceau pervers pur de poids $a$ sur $U$.
Alors $j_{!*}K$ est pervers pur de poids $a$ sur $X$ par [BBD] 5.4.3, 
donc, par [BBD] 5.1.14 et 5.4.1, pour tout $k\in\Z$, 
$\Hp^ki^*j_{!*}K$ est de poids $\leq a+k$
et $\Hp^ki^!j_{!*}K$ est de poids $\geq a+k$.
D'après le lemme \ref{fp_lemme1} de la section 5 à la fin de cette partie, 
$\Hp^ki^*j_{!*}K=0$ si $k\geq 0$ et $\Hp^ki^!j_{!*}K=0$ si $k\leq 0$,
donc, pour tout $k\in\Z$,
$\Hp^ki^*j_{!*}K$ est de poids $\leq a-1$ et $\Hp^ik^!j_{!*}K$ est de poids $\geq a+1$.
Autrement dit, $i^*j_{!*}K\in\DP^{\leq a-1}\subset\DP^{\leq a}$ et $i^!j_{!*}K\in\DP^{\geq a+1}$.

\end{itemize}
\end{proof}

\subsection{Preuve des propositions de la section 3.1}
\hspace{.5cm}

Dans cette section, nous allons prouver les propositions \ref{t_structures:prop1} et \ref{t_structures:prop2} de la section précédente.

\begin{proof}[Démonstration de la proposition \ref{t_structures:prop1}.]
\begin{itemize}
\item[(i)] Il est clair que $\DP^{\leq a}$ et $\DP^{\geq a}$ sont stables par décalage
et que la dualité de Poincaré échange $\DP^{\leq a}$ et $\DP^{\geq -a}$.\newline
Pour montrer la stabilité par extensions de $\DP^{\leq a}$ (resp. $\DP^{\geq a})$,
il suffit de prouver que la catégorie des faisceaux pervers de poids $\leq a$ (resp. $\geq a$)
est une sous-catégorie épaisse de la catégorie des faisceaux pervers,
c'est-à-dire stable par noyaux, conoyaux et extensions.
Par dualité, il suffit de traiter le premier cas.
La stabilité par noyaux et conoyaux résulte de [BBD] 5.3.1,
et la stabilité par extensions se prouve facilement à partir de [BBD] 5.1.9.\newline
Supposons que $a<a'$.
En appliquant la fin de [BBD] 5.1.8 aux objets de cohomologie perverse, on voit que $\DP^{\leq a}\cap\DP^{\geq a'}=0$.
\item[(ii)] Évident.
\item[(iii)] Voir le premier des lemmes ci-dessous. 
\item[(iv)] La condition (ii) de [BBD] 1.3.1 est évidente. 
La condition (i) résulte du point (iii) ci-dessus, et la condition
(iii) est prouvée dans le deuxième des lemmes ci-dessous.
\end{itemize}
\end{proof}

\begin{lemme}\label{lemme_annulation} Soient $X$ un schéma séparé de type fini sur $\Fi_q$, $K,L\in D^b_m(X,\Q_\ell)$ et $a\in\Z$. On suppose que pour tout $i\in\Z$, $\Hp^i(K)$ est de poids $\leq a$ et $\Hp^i(L)$ de poids $\geq a+1$. Alors $R\Hom(K,L)=0.$

\end{lemme}

\begin{proof} Pour tout $i\in\Z$, on a un triangle distingué
$${}^{p}\tau_{\leq i-1}L\fl {}^{p}\tau_{\leq i}L\fl\Hp^i L[-i]\stackrel{+1}{\fl},$$
d'où un triangle distingué
$$R\Hom(K,{}^{p}\tau_{\leq i-1}L)\fl R\Hom(K,{}^{p}\tau_{\leq i}L)\fl R\Hom(K,\Hp^i L)[-i]\stackrel{+1}{\fl}.$$
Si on montre le résultat pour $L$ pervers, on pourra, grâce à ces triangles, en déduire le résultat pour $L$ quelconque en faisant une récurrence sur le cardinal de $\{i\in\Z\mbox{ tq }\Hp^i L\not=0\}$. 
On peut donc supposer $L$ pervers. On se ramène de même au cas où $K$ est pervers.

Notons $F$ le morphisme de Frobenius géométrique.
D'après [BBD] 5.1.2.5, on a pour tout $i\in\Z$ une suite exacte
$$0\fl \Ext^{i-1}(K_{\overline{\Fi}_q},L_{\overline{\Fi}_q})_F\fl \Ext^i(K,L)\fl \Ext^i(K_{\overline{\Fi}_q},L_{\overline{\Fi}_q})^F\fl 0.$$

$K$ est de poids $\leq a$ et $L$ de poids $\geq a+1$, donc $\Ext^i(K_{\overline{\Fi}_q},L_{\overline{\Fi}_q})$ est de poids $>i$ pour tout $i\in\Z$ ([BBD] 5.1.15 (i)). On en déduit que si $i\geq 0$, $\Ext^i(K_{\overline{\Fi}_q},L_{\overline{\Fi}_q})$ est de poids $>0$, donc que
$$\Ext^i(K_{\overline{\Fi}_q},L_{\overline{\Fi}_q})_F=\Ext^i(K_{\overline{\Fi}_q},L_{\overline{\Fi}_q})^F=0$$
(comme dans la preuve de [BBD] 5.1.15). 

D'autre part, $K$ et $L$ sont pervers (donc $K_{\overline{\Fi}_q}$ et $L_{\overline{\Fi}_q}$ aussi), d'où $\Ext^i(K_{\overline{\Fi}_q},L_{\overline{\Fi}_q})=0$ pour $i<0$. 

Finalement, on a obtenu : pour tout $i\in\Z$,
$$\Ext^i(K_{\overline{\Fi}_q},L_{\overline{\Fi}_q})_F=\Ext^i(K_{\overline{\Fi}_q},L_{\overline{\Fi}_q})^F=0.$$
Les suites exactes ci-dessus impliquent que, pour tout $i\in\Z$,
$$\Ext^i(K,L)=0,$$
ce qui est le résultat cherché.

\end{proof}

\begin{remarque} Le lemme ci-dessus reste valable,
avec la même preuve, 
pour une perversité arbitraire (vérifiant les conditions de [BBD] 2.2.1).

\end{remarque}

\begin{lemme}\label{lemme:filtration_poids} Soit $X$ un schéma séparé de type fini sur $\Fi_q$.
Alors pour tout $a\in\Z$ et tout $K\in D^b_m(X,\Q_\ell)$, il existe un triangle distingué dans $D^b_m(X,\Q_\ell)$
$$K_1\fl K\fl K_2\stackrel{+1}{\fl}$$
tel que pour tout $i\in\Z$, $\Hp^i K_1$ soit de poids $\leq a$ et $\Hp^i K_2$ de poids $\geq a+1$.

\end{lemme}

\begin{proof} On fixe $a\in\Z$ et on raisonne par récurrence sur $\card(\{i\in\Z\mbox{ tq }\Hp^iK\not=0\})$.

Supposons qu'il existe $i\in\Z$ tel que $K\simeq\Hp^iK[-i]$.
Alors, d'après [BBD] 5.3.5, il existe un sous-faisceau pervers $L\subset\Hp^iK$ de poids $\leq a$ tel que $\Hp^iK/L$ soit de poids $\geq a+1$.
Il suffit de poser $K_1=L[-i]$ et $K_2=(\Hp^iK/L)[-i]$.

Soit $K\in D^b_m(X,\Q_\ell)$ tel que $n=\card\{i\in\Z\mbox{ tq }\Hp^iK\not=0\}\geq 2$. 
Supposons le lemme prouvé pour tous les $K'\in D^b_m(X,\Q_\ell)$ tels que $\card\{i\in\Z\mbox{ tq }\Hp^iK'\not=0\}<n$, et montrons-le pour $K$.
Il suffit de montrer le résultat suivant : si on a un triangle distingué de $D^b_m(X,\Q_\ell)$
$$K'\fl K\fl K''\stackrel{+1}{\fl}$$
et que la conclusion du lemme vaut pour $K'$ et $K''$, alors elle vaut aussi pour $K$.

On se donne donc des triangles distingués dans $D^b_m(X,\Q_\ell)$
$$\xymatrix{ & & & \\
K'_2\ar[u]^{+1} & & K''_2\ar[u]^{+1} & \\
K'\ar[r]\ar[u] & K\ar[r] & K''\ar[r]^{+1}\ar[u] & \\
K_1'\ar[u] &  & K''_1\ar[u] & }$$
et on suppose que pour tout $i\in\Z$, $\Hp^iK'_1$ et $\Hp^iK''_1$ sont de poids $\leq a$ et $\Hp^iK'_2$ et $\Hp^iK''_2$ de poids $\geq a+1$;
autrement dit, $K'_1$ et $K''_1$ sont dans $\DP^{\leq a}$ et $K'_2$ et $K''_2$ sont dans $\DP^{\geq a+1}$.
D'après le lemme \ref{lemme_annulation}, $R\Hom(K''_1,K'_2[1])=0$, donc il existe un unique morphisme $K''_1\fl K'_1[1]$ qui fait commuter le carré
$$\xymatrix{K''_1\ar[r]\ar[d] & K'_1[1]\ar[d] \\
K''\ar[r] & K'[1]}$$
D'après [BBD] 1.1.11, on peut compléter le diagramme commutatif en traits pleins
pour obtenir un diagramme dont les lignes et les colonnes sont des triangles distingués
et dont tous les carrés sont commutatifs, sauf le carré marqué $-$, qui est anticommutatif :
$$\xymatrix{K''_1[1]\ar[r] & K'_1[2]\ar@{-->}[r] & K_1[2]\ar@{-->}[r]\ar@{}[rd]|{-} & K''_1[2] \\
K''_2\ar@{-->}[r]\ar[u] & K'_2[1]\ar@{-->}[r]\ar[u] & K_2[1]\ar@{-->}[r]\ar@{-->}[u] & K''_2[1]\ar[u] \\
K''\ar[r]\ar[u] & K'[1]\ar[r]\ar[u] & K[1]\ar[r]\ar@{-->}[u] & K''[1]\ar[u] \\
K''_1\ar[u]\ar[r] & K'_1[1]\ar@{-->}[r]\ar[u] & K_1[1]\ar@{-->}[r]\ar@{-->}[u] & K''_1[1]\ar[u]}$$
Comme $\DP^{\leq a}$ et $\DP^{\geq a+1}$ sont stables par extensions (proposition \ref{t_structures:prop1} (i)),
$K_1\in\DP^{\leq a}$ et $K_2\in\DP^{\geq a+1}$.

\end{proof}

\begin{proof}[Démonstration de la proposition \ref{t_structures:prop2}.]
\begin{itemize}
\item[(i)] Soit $K\in D^b_m(X,\Q_\ell)$.
On a un triangle distingué
$$(w_{\leq a+2}K)(1)\fl K(1)\fl (w_{>a+2}K)(1)\stackrel{+1}{\fl}$$
avec $(w_{\leq a+2}K)(1)\in\DP^{\leq a}$ et $(w_{>a+2}K)(1)\in\DP^{>a}$, d'où des isomorphismes canoniques
$w_{\leq a}(K(1))=(w_{\leq a+2}K)(1)$ et $w_{>a}(K(1))=(w_{>a+2}K)(1)$.
\item[(ii)] Soient $K\in D^b_m(X,\Q_\ell)$ et $i\in\Z$.
Comme $\Hp^iw_{\geq a+1}K$ est de poids $\geq a+1$ et que $\Hp^{i+1}w_{\leq a}K$ est de poids $\leq a$,
la flèche $\Hp^iw_{\geq a+1}K\fl \Hp^{i+1}w_{\leq a}K$ est nulle par [BBD] 5.3.1.
\item[(iii) et (iv)] Il suffit de prouver les assertions pour $w_{\leq a}$, celles pour $w_{\geq a}$ en résultant par dualité.\newline
$w_{\leq a}$ commute au décalage parce que $\DP^{\leq a}$ est invariante par décalage.\newline
Soit $K$ un faisceau pervers mixte. 
Si $L$ est le plus grand sous-faisceau pervers de $K$ de poids $\leq a$, alors $K/L$ est de poids $\geq a+1$ ([BBD] 5.3.5).
Donc $w_{\leq a}K=L$, et $w_{\leq a}K$ est pervers.\newline
L'exactitude de la restriction de $w_{\leq a}$ à la catégorie des faisceaux pervers mixtes provient de la fin de [BBD] 5.3.5 
(le fait que les morphismes sont strictement compatibles à la filtration par le poids).\newline
Soient $K\in D^b_m(X,\Q_\ell)$ et $i\in\Z$.
D'après (ii), on a une suite exacte de faisceaux pervers
$$0\fl\Hp^iw_{\leq a}K\fl\Hp^i K\fl\Hp^iw_{\geq a+1}K\fl 0$$
avec $\Hp^iw_{\leq a}K$ de poids $\leq a$ et $\Hp^iw_{\geq a+1}K$ de poids $\geq a+1$,
donc $\Hp^iw_{\leq a}K=w_{\leq a}(\Hp^iK)$.\newline
Le fait que $w_{\leq a}$ envoie les triangles distingués sur des triangles distingués
résulte de la preuve du lemme \ref{lemme:filtration_poids}.
\item[(v)] Les inclusions résultent de [BBD] 4.2.4 et 5.1.14.\newline
Traitons par exemple le cas de $Rf_!$.
Soit $K\in \DP^{\leq a}(X)$.
Pour tous $i,j\in\Z$, $\Hp^jK$ est de poids $\leq a$ par définition de $\DP^{\leq a}(X)$,
donc $Rf_!\Hp^jK$ est de poids $\leq a$ par [BBD] 5.1.14,
et $\Hp^i(Rf_!\Hp^jK)$ est de poids $\leq a+i$ par [BBD] 5.4.1.
Or, par [BBD] 4.2.4, $\Hp^i(Rf_!\Hp^jK)=0$ si $i>d$,
donc $\Hp^i(Rf_!\Hp^jK)$ est de poids $\leq a+d$ pour tous $i,j\in\Z$.
En utilisant la suite spectrale
$$E_2^{ij}=\Hp^i(Rf_!\Hp^jK)\Longrightarrow \Hp^{i+j}Rf_!K,$$
on voit que $\Hp^kRf_!K$ est de poids $\leq a+d$ pour tout $k\in\Z$.

\end{itemize}
\end{proof}

\subsection{t-structures recollées}
\hspace{.5cm}

Nous utiliserons la notion suivante de stratification :

\begin{definition}\label{definition_stratification} Soit $X$ un schéma séparé de type fini sur $\Fi_q$. 
Une \emph{stratification} de $X$ est une partition finie $(S_i)_{0\leq i\leq n}$ de $X$ par des sous-schémas localement fermés (les strates) 
telle que
pour tout $i\in\{0,\dots,n\}$, $S_{i}$ est ouvert dans
$\displaystyle{X-\bigcup_{0\leq j<i}S_{j}}$.

\end{definition}

Soit $X$ un schéma muni d'une stratification $(S_k)_{0\leq k\leq n}$.
Pour tout $k\in\{0,\dots,n\}$, on note $i_k$ l'inclusion de $S_k$ dans $X$.
$U=S_0$ est un ouvert de $X$; on note $j=i_0:U\fl X$ l'inclusion.

La proposition suivante est un cas particulier de [BBD] 1.4.10.

\begin{proposition} Soit $\as=(a_0,\dots,a_n)\in (\Z\cup\{\pm\infty\})^{n+1}$.
On note $\DP^{\leq\as}(X)$ ou $\DP^{\leq\as}$ (resp. $\DP^{>\as}(X)$ ou $\DP^{>\as}$) la sous-catégorie pleine de $D^b_m(X,\Q_\ell)$ 
dont les objets sont les complexes mixtes $K$ tels que
pour tout $k\in\{0,\dots,n\}$,
$i_k^*K\in\DP^{\leq a_k}(S_k)$ (resp. $i_k^!K\in\DP^{>a_k}(S_k)$).

Alors $(\DP^{\leq\as},\DP^{>\as})$ est une t-structure sur $D^b_m(X,\Q_\ell)$.

\end{proposition}

On note aussi
$\DP^{\geq (a_0,\dots,a_n)}=\DP^{>(a_0-1,\dots,a_n-1)}$.
Il est évident d'après la définition de $\DP^{\leq\as}$ et $\DP^{>\as}$ que ce sont des sous-catégories stables par décalage et extensions de $D^b_m(X,\Q_\ell)$,
que la dualité de Poincaré échange $\DP^{\leq\as}$ et $\DP^{\geq -\as}$,
et que 
$$\DP^{\leq (a_0,\dots,a_n)}(1)=\DP^{\leq (a_0-2,\dots,a_n-2)}$$
$$\DP^{>(a_0,\dots,a_n)}(1)=\DP^{>(a_0-2,\dots,a_n-2)}.$$

D'après [BBD] 1.3.3, l'inclusion $\DP^{\leq\as}\subset D^b_m(X,\Q_\ell)$ (resp. $\DP^{>\as}\subset D^b_m(X,\Q_\ell)$) 
admet un adjoint à droite (resp. à gauche), qu'on note
$w_{\leq\as}$ (resp. $w_{>\as}$).
Pour tout $K\in D^b_m(X,\Q_\ell)$, il existe un unique morphisme $w_{>\as}K\fl (w_{\leq\as}K)[1]$ tel que le triangle
$$w_{\leq\as}K\fl K\fl w_{>\as}K\fl (w_{\leq\as}K)[1]$$
soit distingué.\newline
De plus, à isomorphisme unique près, il existe un unique triangle distingué $K'\fl K\fl K''\stackrel{+1}{\fl}$
avec $K'\in\DP^{\leq\as}$ et $K''\in\DP^{>\as}$.

Enfin, la dualité de Poincaré échange $w_{\leq\as}$ et $w_{\geq -\as}$ (car elle échange $\DP^{\leq\as}$ et $\DP^{\geq -\as}$).

\begin{lemme} Si $a_0=\dots=a_n=a$, 
alors $(\DP^{\leq\as},\DP^{>\as})$ est la t-structure $(\DP^{\leq a}(X),\DP^{>a}(X))$ de la section 3.1.

\end{lemme}

\begin{proof} Soit $K\in\DP^{\leq a}$.
D'après le (iv) de la proposition \ref{t_structures:prop2},
pour tout $k\in\{0,\dots,n\}$, $i_k^*K\in\DP^{\leq a}(S_k)$.
Donc $K\in\DP^{\leq\as}$.
Par dualité, on a $\DP^{>a}\subset\DP^{>\as}$.

Soit $K\in\DP^{\leq\as}$.
On a un triangle distingué
$$w_{\leq a}K\fl K\fl w_{>a}K\stackrel{+1}{\fl}.$$
D'après ce qui précède, 
$w_{\leq a}K\in\DP^{\leq\as}$ et $w_{>a}\in\DP^{>\as}$,
donc $w_{\leq a}K=w_{\leq\as}K=K$, et $K\in\DP^{\leq a}$.
Par dualité, on a $\DP^{>\as}\subset\DP^{>a}$.

\end{proof}

\begin{proposition}\label{prop1:t_structures_recollees} Pour tous $a\in\Z\cup\{\pm\infty\}$ et $k\in\{0,\dots,n\}$, on note
$$w_{\leq a}^k=w_{\leq (+\infty,\dots,+\infty,a,+\infty,\dots,+\infty)}$$
$$w_{\geq a}^k=w_{\geq (-\infty,\dots,-\infty,a,-\infty,\dots,-\infty)}$$
où, dans les deux formules, le $a$ est en $k$-ième position (et on commence à compter à $0$).

\begin{itemize}
\item[(i)] On a
$$w_{\leq\as}=w^n_{\leq a_n}\circ\dots\circ w^0_{\leq a_0}$$
$$w_{\geq\as}=w^n_{\geq a_n}\circ\dots\circ w^0_{\geq a_0}.$$
\item[(ii)] Soient $a\in\Z\cup\{\pm\infty\}$, $k\in\{0,\dots,n\}$ et $K\in D^b_m(X,\Q_\ell)$.
Alors on a des triangles distingués (uniques à isomorphisme unique près)
$$w_{\leq a}^kK\fl K\fl Ri_{k*}w_{>a}i_k^*K\stackrel{+1}{\fl}$$
$$i_{k!}w_{\leq a}i_k^!K\fl K\fl w^k_{>a}K\stackrel{+1}{\fl}.$$
\end{itemize}
\end{proposition}

\begin{proof}
\begin{itemize}
\item[(i)] Il suffit d'appliquer plusieurs fois [BBD] 1.4.13.1.
\item[(ii)] Il suffit de traiter le cas de $w_{\leq a}^k$ (celui de $w_{>a}^k$
en résulte par dualité).
On définit $\as\in (\Z\cup\{\pm\infty\})^{n+1}$ par :
$a_r=+\infty$ et si $r\not=k$, et $a_k=a$.

On note
$L=Ri_{k*}w_{>a}i_k^*K\in\DP^{>\as}$,
et on complète le morphisme évident $K\stackrel{adj}{\fl}Ri_{k*}i_k^*K \fl L$ en un triangle distingué
$L'\fl K\fl L\stackrel{+1}{\fl}$.

Il suffit de montrer que $L\in \DP^{>\as}$ et que $L'\in\DP^{\leq\as}$.
Si $r\not=k$, $i_r^!L=0$, car $i_r^!Ri_{k*}=0$,
donc $i_r^!L\in\DP^{>a_r}(S_r)$ (et il est évident que $i_r^*L'\in\DP^{\leq a_r}(S_r)$).
De plus, $i_k^*L=i_k^!L=w_{>a}i_k^*K\in\DP^{>a_k}(S_k)$ et 
on a un triangle distingué
$i_k^*L'\fl i_k^*K\fl i_k^*L=w_{>a}i_k^*K\stackrel{+1}{\fl}$,
donc $i_k^*L'\simeq w_{\leq a}i_k^*K\in\DP^{\leq a_k}(S_k)$.

\end{itemize}
\end{proof}

\begin{theoreme}\label{th:simplification_formule_traces} 
Pour tout $\as\in (\Z\cup\{\pm\infty\})^{n+1}$,
pour tout $K\in\DP^{\leq a_0}(U)$ 
on a
$$[w_{\leq\as}Rj_*K]=\sum_{1\leq n_1<\dots <n_r\leq n}(-1)^{r}[Ri_{n_r*}w_{>a_{n_r}}i_{n_r}^*\dots Ri_{n_1*}w_{>a_{n_1}}i_{n_1}^*Rj_*K]$$
dans le groupe de Grothendieck de $D^b_m(X,\Q_\ell)$.

\end{theoreme}

\begin{proof} 

D'après la proposition ci-dessus, on a, dans l'anneau des endomorphismes du groupe 
de Grothendieck de $D^b_m(X,\Q_\ell)$ :
$$w_{\leq\as}=w_{\leq a_n}^n\circ\dots\circ w_{\leq a_0}^0=(1-Ri_{n*}w_{>a_n}i_n^*)\circ\dots\circ (1-Ri_{1*}w_{>a_1}i_1^*)\circ(1-Rj_*w_{>a_0}j^*).$$
Le théorème résulte de cette égalité et du fait que $Rj_*w_{>a_0}j^*Rj_*K=0$ 
(car $K\in\DP^{\leq a_0}(U)$).

\end{proof}

\subsection{Propriétés supplémentaires des t-structures recollées}
\hspace{.5cm}

\begin{proposition}\label{prop2:t_structures_recollees}
\begin{itemize}
\item[(i)] Si $K\in D^{\leq\as}$ et $L\in\DP^{>\as}$, alors $R\Hom(K,L)=0$.
\item[(ii)] $w_{\leq\as}$ et $w_{>\as}$ commutent au foncteur de décalage $[1]$,
et ils envoient les triangles distingués sur des triangles distingués.\newline
Pour tout $K\in D^b_m(X,\Q_\ell)$, on a
$$w_{\leq (a_0,\dots,a_n)}(K(1))=(w_{\leq (a_0+2,\dots,a_n+2)}K)(1)$$
$$w_{>(a_0,\dots,a_n)}(K(1))=(w_{>(a_0+2,\dots,a_n+2)}K)(1).$$
\item[(iii)] Soit $\as'=(a'_0,\dots,a'_n)\in (\Z\cup\{\pm\infty\})^{n+1}$ tel que $a_k\leq a'_k$
pour tout $k\in\{1,\dots,n\}$.
Alors, pour tous $K\in\DP^{\leq\as}$ et $L\in\DP^{>\as'}$, le morphisme canonique
$$R\Hom(K,L)\fl R\Hom(j^*K,j^*L)$$
est un isomorphisme.
\item[(iv)] Soient $f:Y\fl X$ un morphisme et $(S'_k)_{0\leq k\leq n}$ une stratification de $Y$
telle que, pour tout $k\in\{0,\dots,n\}$, $f(S'_k)\subset S_k$.
On suppose que la dimension des fibres de $f$ est inférieure ou égale à $d$.
Alors
$$f^*(\DP^{\leq (a_0,\dots,a_n)}(X))\subset\DP^{\leq (a_0+d,\dots,a_n+d)}(Y)$$
$$f^!(\DP^{>(a_0,\dots,a_n)}(X))\subset\DP^{>(a_0-d,\dots,a_n-d)}(Y)$$
$$Rf_*(\DP^{>(a_0,\dots,a_n)}(Y))\subset\DP^{>(a_0-d,\dots,a_n-d)}(X)$$
$$Rf_!(\DP^{\leq (a_0,\dots,a_n)}(Y))\subset\DP^{\leq (a_0+d,\dots,a_n+d)}(X).$$
\end{itemize}
\end{proposition}

\begin{proof}
\begin{itemize}
\item[(i)] Montrons le résultat par récurrence sur $n$.
Si $n=0$, c'est le lemme \ref{lemme_annulation}.
Soit $n\geq 1$, et supposons le résultat vrai pour $n'<n$.
On note $\as'=(a_0,\dots,a_{n-1})$, $V=\displaystyle{\bigcup_{k=0}^{n-1}S_k}$, $Y=S_n=X-V$, 
$j_1$ l'immersion ouverte de $V$ dans $X$ et 
$i$ l'immersion fermée de $S_n$ dans $X$.
Soient $K\in\DP^{\leq\as}$ et $L\in\DP^{>\as}$.
On a un triangle distingué
$$R\Hom(i^*K,i^!L)\fl R\Hom(K,L)\fl R\Hom(j_1^*K,j_1^*L)\stackrel{+1}{\fl}.$$
Or $j_1^*K\in\DP^{\leq\as'}(U)$, $j_1^*L\in\DP^{>\as'}(U)$, $i^*K\in\DP^{\leq a_n}(Y)$ et $i_k^!L\in\DP^{>a_n}(Y)$,
donc, d'après l'hypothèse de récurrence,
$$R\Hom(j_1^*K,j_1^*L)=R\Hom(i^*K,i^!L)=0,$$
d'où 
$$R\Hom(K,L)=0.$$
\item[(ii)] $w_{\leq\as}$ et $w_{>\as}$ commutent au foncteur de décalage car $\DP^{\leq\as}$ et $\DP^{>\as}$ sont stables par décalage.
Soit $K\fl K'\fl K''\stackrel{+1}{\fl}$ un triangle distingué.
D'après [BBD] 1.1.11, on peut construire un diagramme commutatif dont les lignes et les colonnes sont distinguées
$$\xymatrix{w_{\leq\as}K\ar[r]\ar[d] & w_{\leq\as}K'\ar[r]\ar[d] & L\ar[r]^{+1}\ar[d] & \\
K\ar[r]\ar[d] & K'\ar[r]\ar[d] & K''\ar[r]^{+1}\ar[d] & \\
w_{>\as}K\ar[r]\ar[d]^{+1} & w_{>\as}K'\ar[r]\ar[d]^{+1} & L'\ar[r]^{+1}\ar[d]^{+1} & \\
 & & & }$$
Comme $\DP^{\leq\as}$ et $\DP^{>\as}$ sont des sous-catégories stables par extensions de $D^b_m(X,\Q_\ell)$,
$L\in\DP^{\leq\as}$ et $L'\in\DP^{>\as}$, donc 
$L=w_{\leq\as}K''$ et $L'=w_{>\as}K''$.\newline
La dernière assertion se prouve exactement comme la propriété analogue dans le (i) de la proposition \ref{t_structures:prop2}.
\item[(iii)] Notons $i$ l'immersion fermée $X-U=S_1\cup\dots\cup S_n\subset X$, $\underline{b}=(a_1,\dots,a_n)$ et $\underline{b}'=(a'_1,\dots,a'_n)$.
Soient $K\in\DP^{\leq\as}$ et $L\in\DP^{>\as'}$.
On a un triangle distingué canonique
$$R\Hom(i^*K,i^!L)\fl R\Hom(K,L)\fl R\Hom(j^*K,j^*L)\stackrel{+1}{\fl}.$$
Or $i^*K\in\DP^{\leq\underline{b}}(X-U)$ et $i^!L\in\DP^{>\underline{b}'}(X-U)$,
donc, d'après le point (i), $R\Hom(i^*K,i^!L)=0$.
\item[(iv)] Il suffit de traiter les cas de $f^*$ et $Rf_!$, car
ceux de $f^!$ et $Rf_*$ en résultent par dualité.\newline
Pour tout $k\in\{0,\dots,n\}$, on note $i'_k$ l'inclusion $S'_k\subset Y$ et 
$f_k:S'_k\fl S_k$ la restriction de $f$.
Soit $K\in\DP^{\leq\as}(X)$.
Pour tout $k\in\{0,\dots,n\}$, ${i'_k}^*f^*K=f_k^*i_k^*K$; 
$i_k^*K\in\DP^{\leq a_k}(S_k)$ par la définition de $\DP^{\leq\as}(X)$,
donc, d'après le (iv) de la proposition \ref{t_structures:prop2}, $f_k^*i_k^*K\in\DP^{\leq a+d}(S'_k)$.
On a donc bien $f^*K\in\DP^{\leq (a_0+d,\dots,a_n+d)}(Y)$.

Soit $K\in\DP^{\leq\as}(Y)$.
Fixons $k\in\{0,\dots,n\}$. 
Comme $S'_k=f^{-1}(S_k)$,
le diagramme suivant est cartésien aux nilpotents près
$$\xymatrix{S'_k\ar[r]^{i'_k}\ar[d]_{f_k} & Y\ar[d]_f \\
S_k\ar[r]^{i_k} & X}$$
donc, d'après le théorème de changement de base propre, ${i_k}^*Rf_!K\simeq Rf_{k!}{i'_k}^*K$.
Or ${i'_k}^*K\in\DP^{\leq a_k}(S'_k)$, donc, d'après le (iv) de la proposition \ref{t_structures:prop2},
$i_k^*Rf_!K\simeq Rf_{k!}{i'_k}^*K\in\DP^{\leq a_k+d}(S_k)$.
On a donc bien $Rf_!K\in\DP^{\leq (a_0+d,\dots,a_n+d)}(X)$.
\end{itemize}
\end{proof}

La proposition suivante est une reformulation de [BBD] 1.4.14 dans le cas particulier considéré.

\begin{proposition}\label{prop3:t_structures_recollees}
Soit $\as=(a_0,\dots,a_n)\in (\Z\cup\{\pm\infty\})^{n+1}$.
On note $\protect{\as'\!=\!(a_0,a_1+1,\dots,a_n+1)}$.
Alors, pour tout $K\in\DP^{\leq a_0}(U)\cap\DP^{\geq a_0}(U)$, 
$w_{\geq\as'}j_!K=w_{\leq\as}Rj_*K$ est l'unique prolongement de $K$ dans $\DP^{\leq\as}\cap\DP^{\geq\as'}$.

En particulier, si $K$ est pervers pur de poids $a$ sur $U$, on a
$$w_{\geq (a,a+1,\dots,a+1)}j_!K=j_{!*}K=w_{\leq a}Rj_*K,$$
et par dualité on obtient aussi
$$w_{\geq a}j_!K=j_{!*}K=w_{\leq (a,a-1,\dots,a-1)}Rj_*K.$$
(On retrouve le résultat du théorème \ref{th:calcul_IC_general}.)

\end{proposition}

\subsection{Quelques lemmes techniques}
\hspace{.5cm}

Ce paragraphe contient quelques lemmes utilisés dans les preuves des résultats des parties 3 et 4.

\begin{lemme}\label{fp_lemme1} Soient $X$ un schéma séparé de type fini sur $\Fi_q$ et $(S_\alpha)$ une partition finie de $X$ par des sous-schémas localement fermés. 
Pour tout $\alpha$, on note $i_\alpha:S_\alpha\fl X$ l'inclusion. Alors, pour tout $K\in D^b_c(X,\Q_\ell)$, on a
\begin{itemize}
\item[(a)] $K\in {}^{p}D_c^{\leq 0}$ si et seulement si pour tout $\alpha$, pour tout $i\geq 1$, on a $\Hp^i(i_\alpha^*K)=0$;
\item[(b)] $K\in {}^{p}D_c^{\geq 0}$ si et seulement si pour tout $\alpha$, pour tout $i\leq -1$, on a $\Hp^i(i_\alpha^!K)=0$.

\end{itemize}
De plus, si $U$ est un ouvert de $X$ réunion de strates, si $j:U\fl X$ est l'inclusion et si $K$ est un faisceau pervers sur $U$, alors $j_{!*}K$ est l'unique prolongement $L\in D^b_c(X,\Q_\ell)$ de $K$ tel que : pour tout $\alpha$ tel que $S_\alpha\subset X-U$, on a
$$\left\{\begin{array}{c}{}^{p}H^i(i_\alpha^*L)=0\mbox{ pour }i\geq 0 \\
{}^{p}H^i(i_\alpha^!L)=0\mbox{ pour }i\leq 0\end{array}\right..$$

\end{lemme}

\begin{proof} Il suffit de montrer (a), car (b) en résulte par dualité. 
D'après [BBD] 2.2.5, les $i_\alpha^*$ sont $t$-exacts à droite, donc, si $K\in {}^{p}D_c^{\leq 0}$, on a bien $\Hp^i(i_\alpha^*K)=0$ pour tout $\alpha$ et pour $i\geq 1$.
 
Réciproquement, soit $K\in D^b_c(X,\Q_\ell)$ tel que pour tout $\alpha$, pour tout $i\geq 1$, $\Hp^i(i_\alpha^*K)=0$.
On sait par [BBD] 2.2.12 qu'un complexe $L$ est dans ${}^{p}D_c^{\leq 0}$ si et seulement si pour tout point $x$ de $X$, notant $i_x:x\fl X$ et $dim(x)=dim(\overline{\{x\}})$, on a $H^i(i_x^*L)=0$ pour $i<p(2 dim(x))$. On va utiliser cette caractérisation pour montrer que $K\in {}^{p}D_c^{\leq 0}$.
Soit $x$ un point de $X$, et soit $\alpha$ tel que $x\in S_\alpha$. Comme $S_\alpha$ est localement fermé dans $X$, $\overline{\{x\}}\cap S_\alpha$ est ouvert dans $\overline{\{x\}}$, donc $dim(\overline{\{x\}}\cap S_\alpha)=dim(\overline{\{x\}})$ ($\overline{\{x\}}$ est irréductible), et $dim(x)$ ne change pas si on considère $x$ comme un point de $S_\alpha$. 
Comme par hypothèse $i_\alpha^*K\in {}^{p}D_c^{\leq 0}(S_\alpha)$, on a bien $H^i(i_x^*K)=0$ si $i<p(2 dim(x))$. 

Montrons enfin la dernière assertion du lemme. $U$ est un ouvert de $X$ réunion de strates, $j:U\fl X$ est l'inclusion, $K$ est un faisceau pervers sur $U$. D'après [BBD] 1.4.24 (qui s'applique par [BBD] 2.2.3 et 2.2.11), on a $\Hp^0(i_\alpha^*j_{!*}K)=0$ pour tout $\alpha$ tel que $S_\alpha\subset X-U$. L'annulation des $\Hp^0(i_\alpha^!j_{!*}K)$ s'en déduit par dualité.

Soit $L\in D^b_c(X,\Q_\ell)$, muni d'un isomorphisme $j^*L\simeq K$, tel que, pour tout $\alpha$ tel que $S_\alpha\subset X-U$, on ait $\Hp^i(i_\alpha^*L)=0$ pour $i\geq 0$ et $\Hp^i(i_\alpha^!L)=0$ pour $i\leq 0$. D'après ce qui précède, on sait que $L$ est pervers.
En raisonnant par récurrence sur le cardinal de $\{\alpha\mbox{ tq }S_\alpha\subset X-U\}$, on se ramène au cas où $X-U=S_\alpha$ est une strate. Notons $i=i_\alpha$. 
On a un triangle distingué
$$i_*i^!L\fl L\fl Rj_*j^*L\simeq Rj_* K\stackrel{+1}{\fl},$$
d'où une suite exacte
$$\Hp^0(i_*i^! L)\fl \Hp^0(L)=L\fl \Hp^0(Rj_*K).$$
Comme $i_*$ est t-exact ([BBD] 2.2.6), $\Hp^0(i_*i^!L)=i_*\Hp^0(i^!L)=0$, et le morphisme $L\fl\Hp^0(Rj_*K)$ est injectif.
D'autre part, on a un triangle distingué
$$j_!j^*L\simeq j_!K\fl L\fl i_*i^*L\stackrel{+1}{\fl},$$
d'où une suite exacte
$$\Hp^0(j_!K)\fl L\fl \Hp^0(i_*i^* L)=i_*\Hp^0(i^* L)=0.$$
Le morphisme $\Hp^0(j_!K)\fl L$ est donc surjectif, ce qui finit la démonstration.

\end{proof}

\begin{lemme}\label{fp_lemme2} Soient $X$ un schéma de type fini lisse purement de dimension $d$ sur un corps $k$ de caractéristique $0$ ou fini et $K\in D^b_c(X,\Q_\ell)$. On suppose que les $H^i(K)$ sont lisses. Alors, pour tout $i\in \Z$, $\Hp^i(K)=H^{i-d}(K)[d]$.

\end{lemme}

\begin{proof} On montre le résultat par récurrence sur le cardinal $N(K)$ de\newline
 $\{i\in\Z\mbox{ tq }H^i(K)\not=0\}$.

Si $N(K)=1$, on a $K\simeq H^i(K)[-i]$ pour un $i\in\Z$, donc $K[i+d]$ est pervers, et
$$\Hp^j(K)=\left\{\begin{array}{ll}0 & \mbox{ si }j\not=i+d \\
H^i(K)[d] & \mbox{ si }j=i+d\end{array}\right..$$

Soit $K$ tel que $N(K)>1$, et supposons le résultat prouvé pour tous les $L$ tels que\newline
 $N(L)<N(K)$. Soit $i=max\{k\in\Z\mbox{ tq }H^k(K)\not=0\}$. Comme $H^i(K)[-i]$ est lisse, on a comme plus haut
$$\Hp^j(H^i(K)[-i])=\left\{\begin{array}{ll}0 & \mbox{ si }j\not=i+d \\
H^i(K)[d] & \mbox{ si }j=i+d\end{array}\right..$$
D'autre part, d'après l'hypothèse de récurrence, on a 
$$\Hp^j(\tau_{\leq i-1}K)=H^{j-d}(\tau_{\leq i-1}K)[d]=\left\{\begin{array}{ll}0 & \mbox{ si }j\geq i+d \\
H^{j-d}(K)[d] & \mbox{ si }j<i+d\end{array}\right..$$
On conclut en utilisant la suite exacte longue de cohomologie perverse du triangle distingué 
$$\tau_{\leq i-1}K\fl\tau_{\leq i}K\simeq K\fl H^i(K)[-i]\stackrel{+1}{\fl}.$$

\end{proof}

\section{Complexes pondérés sur les compactifications de Baily-Borel}

Dans la suite, $(\G,\X)$ est une des données de Shimura $(\GSp_{2d,\Q},\X_d)$ de la section 1.1,
$\K=\K_d(n)$ avec $n\geq 3$,
et $\ell$ est un nombre premier.
On choisit un nombre premier $p\not=\ell$ qui ne divise pas $n$ et 
on travaille sur les réductions modulo $p$ des schémas de la partie 1,
qu'on notera $M^{\K}(\G,\X)$, $M^{\K}(\G,\X)^*$, etc.
Pour tout $V\in D^b(Rep_{\G(\Q_\ell)})$, on notera encore $\F^{\K}V$
le complexe $\ell$-adique mixte sur $M^{\K}(\G,\X)$ obtenu en réduisant le complexe
$\F^{\K}V$ sur $\Mod_{d,n}[1/\ell]$.

\subsection{Complexes pondérés}
\hspace{.5cm}

Dans cette section, on définit à l'aide des foncteurs $w_{\leq\as}$ de 3.3
une famille de complexes sur la compactification de Baily-Borel,
qu'on appellera complexes pondérés,
et dont les complexes d'intersection sont des cas particuliers.
Ces complexes pondérés sont des analogues en caractéristique finie 
des complexes pondérés définis sur $M^{\K}(\G,\X)^*(\C)$ par
Goresky, Harder et MacPherson dans [GHM].

\begin{notation} Soit $t\in\Z$.
Pour tout $V\in Rep_{\G(\Q_\ell)}$, on note $w_{< t}V$ (resp. $w_{\geq t}V$) la plus grande sous-représentation de $V$
de poids $<t$ (resp. $\geq t$).

Les foncteurs exacts $w_{< t}$ et $w_{\geq t}$ s'étendent trivialement à la catégorie dérivée $D^b(Rep_{\G(\Q_\ell)})$.
En effet, comme $Rep_{\G(\Q_\ell)}$ est une catégorie semi-simple
(car $\G$ est réductif), le foncteur de cohomologie $H^*$ est une équivalence de catégories de $D^b(Rep_{\G(\Q_\ell)})$
avec la catégorie des objets gradués de $Rep_{\G(\Q_\ell)}$ qui sont nuls en degré
assez grand et en degré assez petit.

Pour tout $V\in D^b(Rep_{\G(\Q_\ell)})$, on a
$$V=w_{< t}V\oplus w_{\geq t}V,$$
et, pour tout $i\in\Z$, $H^i(w_{< t}V)=w_{< t}H^i(V)$ est de poids $< t$ et $H^i(w_{\geq t}V)=w_{\geq t}H^i(V)$ de poids $\geq t$.

\end{notation}

\begin{lemme}\label{relation_tronques} On note $c=d(d+1)/2$ (c'est la dimension de $M^{\K}(\G,\X)$).
Alors, pour tous $V\in D^b(Rep_{\G(\Q_\ell)})$ et $a\in\Z\cup\{\pm\infty\}$,
$$w_{\leq a}\F^{\K}V=\F^{\K}w_{\geq c-a}V$$
$$w_{>a}\F^{\K}V=\F^{\K}w_{<c-a}V.$$

\end{lemme}

\begin{proof} Notons $t=c-a$, $K=\F^{\K}V$, $K_1=\F^{\K}w_{\geq t}V$ et $K_2=\F^{\K}w_{<t}V$.
Comme $V=w_{<t}V\oplus w_{\geq t}V$, on a $K=K_1\oplus K_2$.
Il suffit de montrer que $K_1$ est dans $\DP^{\leq a}(M^{\K}(\G,\X))$ et $K_2$ dans $\DP^{>a}(M^{\K}(\G,\X))$.

$M^{\K}(\G,\X)$ est lisse et pour tout $i\in\Z$, $H^i(K_1)=\F^{\K}H^i(w_{\geq t}V)$ et $H^i(K_2)=\F^{\K}H^i(w_{<t}V)$ sont lisses,
donc, d'après le lemme \ref{fp_lemme2}, pour tout $i\in\Z$, $\Hp^i(K_1)=H^{i-c}(K_1)[c]$ et $\Hp^i(K_2)=H^{i-c}(K_2)[c]$.

Soit $i\in\Z$.
D'après la proposition \ref{poids_faisceaux}, $H^{i-c}(K_1)$ est de poids $\leq -t$ et $H^{i-c}(K_2)$ est de poids $>-t$.
On en déduit que $\Hp^i(K_1)=H^{i-c}(K_1)[c]$ est de poids $\leq -t+c=a$ et que $\Hp^i(K_2)=H^{i-c}(K_2)[c]$ est de poids $>-t+c=a$.

\end{proof}

On pose $M^*=M^{\K}(\G,\X)^*$, $M_0=M^{\K}(\G,\X)$ et, pour tout $r\in\{1,\dots,d\}$, on note $M_r$ l'union des strates de bord correspondant à $(\QP_{d-r},\Y_{d-r})$.
Pour tout $r\in\{0,\dots,d\}$, 
on pose $c_r=dim(M_r)=(d-r)(d+1-r)/2$.

$(M_0,M_1,\dots,M_q)$ est une stratification de $M^*$ (au sens de la définition \ref{definition_stratification})
, et c'est toujours celle qu'on utilisera dans la suite.

\begin{definition}\label{def:complexes_ponderes} Soient $t_0,\dots,t_{d-1}\in\Z\cup\{\pm\infty\}$.
Pour tout $r\in\{1,\dots,d\}$, on pose $a_r=-t_{d-r}+c_r$.
On définit un foncteur additif triangulé
$$W^{\geq t_0,\dots,\geq t_{d-1}}:D^b(Rep_{\G(\Q_\ell)})\fl D^b_m(M^{\K}(\G,\X)^*)$$
de la manière suivante : pour tout $m\in\Z$, si $V\in D^b(Rep_{\G(\Q_\ell)})$ est tel que $H^i(V)$ soit de poids $m$ pour tout $i\in\Z$, alors
$$W^{\geq t_0,\dots,\geq t_{d-1}}V=w_{\leq (-m+c_0,-m+a_1,\dots,-m+a_d)}Rj_*\F^{\K}V.$$

\end{definition}

\begin{definition}\label{def:IC} Soit $V\in Rep_{\G(\Q_\ell)}$.
Comme $M^{\K}(\G,\X)$ est de dimension $c_0$, $\F^{\K}V[c_0]$ est un faisceau pervers sur $M^{\K}(\G,\X)$.
On pose
$$IC^{\K}V=(j_{!*}(\F^{\K}V[c_0]))[-c_0].$$

\end{definition}

\begin{proposition}\label{calcul_IC} 
\begin{itemize}
\item[(1)] Pour tous $t_0,\dots,t_{d-1}\in\Z\cup\{\pm\infty\}$ et pour tout $V\in D^b(Rep_{\G(\Q_\ell)})$, on a un isomorphisme canonique
$$D(W^{\geq t_0,\dots,\geq t_{d-1}}V)\simeq W^{\geq s_0,\dots,\geq s_{d-1}}(V^*)[2c_0](c_0),$$
où $V^*=R\Hom(V,\Q_\ell)$ est la représentation duale de $V$ et $s_r=1-t_r+2(c_{d-r}-c_0)$ pour tout $r\in\{0,\dots,d-1\}$.
\item[(2)] Notons, pour tout $r\in\{0,\dots,d-1\}$, $t_r=1+c_{d-r}-c_0=1-codim(M_{d-r})$ et $s_r=c_{d-r}-c_0=-codim(M_{d-r})$.
Alors, pour tout $V\in Rep_{\G(\Q_\ell)}$, on a des isomorphismes canoniques
$$IC^{\K}V\simeq W^{\geq t_0,\dots,\geq t_{d-1}}V\iso W^{\geq s_0,\dots,\geq s_{d-1}}V.$$

\end{itemize}
\end{proposition}

\begin{proof}
\begin{itemize}
\item[(1)] On peut supposer que $V\in Rep_{\G(\Q_\ell)}$ et que $V$ est pure de poids $m\in\Z$.
$V^*$ est alors une représentation de $\G(\Q_\ell)$ pure de poids $-m$.
On pose $a_0=-m+c_0$ et, pour tout $r\in\{1,\dots,d\}$, $a_r=-(t_{d-r}+m)+c_r$.
Alors
$$W^{\geq t_0,\dots,\geq t_{d-1}}V=w_{\leq (a_0,\dots,a_d)}Rj_*\F^{\K}V,$$
donc
$$\begin{array}{rcl}D(W^{\geq t_0,\dots,\geq t_{d-1}}V) & = & w_{\geq (-a_0,\dots,-a_d)}(j_!\F^{\K}V^*[2c_0](c_0)) \\
& = & (w_{\geq (-a_0+2c_0,\dots,-a_d+2c_0)}j_!\F^{\K}V^*)[2c_0](c_0).\end{array}$$
D'après la proposition \ref{prop3:t_structures_recollees},
$$w_{\geq (-a_0+2c_0,\dots,-a_d+2c_0)}j_!\F^{\K}V^*=w_{\leq (-a_0+2c_0,-a_1+2c_0-1,\dots,-a_d+2c_0-1)}Rj_*\F^{\K}V^*.$$
Notons $b_0=m+c_0$ et, pour tout $r\in\{1,\dots,d\}$, $b_r=-(s_{d-r}-m)+c_r$.
Alors $b_0=-a_0+2c_0$ et, pour tout $r\in\{1,\dots,d\}$, $b_r=-a_r-1+2c_0$.
Donc
$$\begin{array}{rcl}D(W^{\geq t_0,\dots,\geq t_{d-1}}V) & = & (w_{\leq (b_0,\dots,b_d)}Rj_*\F^{\K}V^*)[2c_0](c_0) \\
& = & W^{\geq s_0,\dots,\geq s_{d-1}}(V^*)[2c_0](c_0).\end{array}$$
\item[(2)] On peut supposer que $V$ est pure. Soit $m$ son poids.
$\F^{\K}V$ est lisse pur de poids $-m$ sur $M_0$, qui est lisse de dimension $c_0$, 
donc le seul faisceau de cohomologie perverse non nul de $\F^{\K}V$
est $\Hp^{c_0}\F^{\K}V=\F^{\K}V[c_0]$, qui est pur de poids $-m+c_0$.
D'après la proposition \ref{prop3:t_structures_recollees},
$$\begin{array}{rcl}IC^{\K}V[c_0] & = & j_{!*}(\F^{\K}V[c_0]) \\
& = & w_{\leq (-m+c_0,\dots,-m+c_0)}Rj_*\F^{\K}V[c_0] \\
& = & w_{\leq (-m+c_0,-m-1+c_0,\dots,-m-1+c_0)}Rj_*\F^{\K}V[c_0],\end{array}$$
donc 
$$\begin{array}{rcl}IC^{\K}V & = & w_{\leq (-m+c_0,\dots,-m+c_0)}Rj_*\F^{\K}V \\
& = & w_{\leq (-m+c_0,-m-1+c_0,\dots,-m-1+c_0)}Rj_*\F^{\K}V.\end{array}$$
Pour conclure, il suffit de remarquer que, pour tout $r\in\{1,\dots,d\}$,\newline
$-t_{d-r}-m+c_r=-m-1+c_0$ et $-s_{d-r}-m+c_r=-m+c_0$.

\end{itemize}
\end{proof}

\subsection{Restrictions des complexes pondérés aux strates}
\hspace{.5cm}

On note $\Se=\Gr_{m,\Q}.I_{2d}$ le centre de $\G$, pour tout $r\in\{1,\dots,d-1\}$, 
$$\Se_r=\left\{\left(\begin{array}{ccc}\lambda^2 I_{d-r} & 0 & 0 \\
0 & \lambda I_{2r} & 0 \\
0 & 0 & I_{d-r}\end{array}\right),\lambda\in\Gr_{m,\Q}\right\}$$
et
$$\Se_0=\left(\begin{array}{cc}\Gr_{m,\Q}.I_d & 0 \\
0 & I_d\end{array}\right).$$
Pour tout $r\in\{0,\dots,d-1\}$, $\Se_r\simeq\Gr_{m,\Q}$ est le centre de $\G_r$.

On utilisera les notations de 1.2 et les notations suivantes :
Soit $S\subset\{0,\dots,d-1\}$ non vide.
On pose
$$\Pa_S=\bigcap_{s\in S}\Pa_s.$$
C'est un sous-groupe parabolique standard de $\G$,
dont on note $\N_S$ le radical unipotent et $\Le_S=\Pa_S/\N_S$ le quotient de Levi.
Soit $r=min(S)$.
On a $\QP_r\subset\Pa_S\subset\Pa_r$, donc
$\Pa_S/\N_r$ est produit direct de $\G_r$ et d'un sous-groupe parabolique $\Pa_{\lin,S}$
de $\Le_{\lin,r}$,
dont on note $\Le_{\lin,S}$ le quotient de Levi.

On pose
$$\Hr_{S}=\K \cap\Pa_S(\Q)\QP_r(\Af)=\K \cap\Pa_{\lin,S}(\Q)\QP_r(\Af)$$
$$\Hr_{\lin,S}=\K\cap\Pa_{\lin,S}(\Q)\N_r(\Af)$$
$$\Gamma_{\ell,S}=\Hr_{\lin,S}/(\K\cap\N_S(\Q)\N_r(\Af))=(\K\cap\Pa_{\lin,S}(\Q)\N_S(\Af))/(\K\cap\N_S(\Af)).$$
On a $\K_r=\Hr_S/\Hr_{\lin,S}$, et
$\Gamma_{\lin,S}=Ker(\Le_{\lin,S}(\Z)\fl\Le_{\lin,S}(\Z/n\Z))$ est un sous-groupe arithmétique net de $\Le_{\lin,S}(\Q)$.

Soit $V\in Rep_{\Le_S(\Q_\ell)}$.
Si $(t_s)_{s\in S}\in (\Z\cup\{\pm\infty\})^S$, on note
$$V_{<t_s,s\in S}$$
le sous-espace vectoriel de $V$ sur lequel, pour tout $s\in S$, $\Se_s$ agit par des caractères $x\fle x^{t}$ avec $t<t_s$.
Comme les $\Se_s$ sont centraux dans $\Le_S$, $V_{<t_s,s\in S}$ est stable par $\Le_S(\Q_\ell)$.

La définition ci-dessus s'étend trivialement aux complexes et donne un foncteur exact
$$\left\{\begin{array}{rcl}D^b(Rep_{\Le_S(\Q_\ell)}) & \fl & D^b(Rep_{\Le_S(\Q_\ell)}) \\
V & \fle & V_{<t_s,s\in S}\end{array}\right..$$ 

On définit de même $V_{\geq t_s,s\in S}$.

\begin{theoreme}\label{th:restriction_W_bord} Soient $t_0,\dots,t_{d-1}\in\Z\cup\{\pm\infty\}$ et $V\in D^b(Rep_{\G(\Q_\ell)})$. On suppose que tous les $H^i(V)$ sont purs de même
poids $m\in\Z$.
On fixe $r\in\{0,\dots,d-1\}$ et $g\in\G(\widehat{\Z})$. 
Alors
\vspace{1cm}
\begin{flushleft}$\displaystyle{[i_{g,r}^*W^{\geq t_0,\dots,\geq t_{d-1}}(V)]=}$\end{flushleft}
\begin{flushright}$\displaystyle{\F^{\K_r}\left(\sum_{S}\sum_{i\in I_S}(-1)^{\card(S)-1}\left[R\Gamma\left(\Gamma_{S},R\Gamma(Lie(\N_{S}),g_ig.V)_{\geq t_{r},< t_{s},s\in S\sous\{r\}}\right)\right]\right)}$\end{flushright}
où $S$ parcourt les sous-ensembles de $\{r,\dots,d-1\}$ contenant $r$ 
et, pour tout
$r\in S\subset\{r,\dots,d-1\}$,
$(g_i)_{i\in I_S}$ est un système de représentants dans $\Pa_r(\Z)\QP_r(\widehat{\Z})$
du double quotient $\Pa_{S}(\Q)\QP_r(\Af)\sous\Pa_r(\Q)\QP_r(\Af)/\Hr_{r}$.

Supposons que $V$ est concentré en degré $0$. Alors,
si $t_s=[s(s+1)-d(d+1)]/2$ pour tout $s\in\{0,\dots,d-1\}$, ou si
$t_s=1+[s(s+1)-d(d+1)]/2$ pour tout $s\in\{0,\dots,d-1\}$,
on a $IC^{\K}V=W^{\geq t_0,\dots,\geq t_{d-1}}(V)$ (proposition \ref{calcul_IC} (2)).
On obtient donc en particulier une formule pour $[i_{g,r}^*IC^{\K}V]$.

\end{theoreme}

\begin{proof} Le théorème résulte du théorème \ref{th:simplification_formule_traces} et
de la proposition ci-dessous.

\end{proof}

\begin{remarque} Grâce à ce théorème et aux résultats de [K],
on peut calculer la fonction trace de Frobenius pour un complexe pondéré
(et en particulier pour les complexes d'intersection à coefficients dans
une représentation $V$ de $\G(\Q_\ell)$).

\end{remarque}

Rappelons qu'on a noté, pour tout $r\in\{1,\dots,d\}$, $M_r$ l'union des strates
$Im(i_{g,d-r})$, $g\in\G(\widehat{\Z})$ (c'est-à-dire l'union des strates de 
$M^{\K}(\G;\X)^*$ associées au sous-groupe parabolique maximal $\Pa_{d-r}$
de $\G$).
On note $i_r$ l'inclusion de $M_r$ dans $M^{\K}(\G,\X)^*$.

\begin{proposition}\label{restr_tronques_bord} Soient $r_1,\dots,r_c\in\{0,\dots,d-1\}$ tels que 
$r_1>\dots>r_c$,
$a_1,\dots,a_c\in\Z\cup\{\pm\infty\}$, $V\in D^b(Rep_{\G(\Q_\ell)})$, $r\in\{0,\dots,r_c\}$ et $g\in\G(\widehat{\Z})$.
Pour tout $i\in\{1,\dots,c\}$, on pose $t_i=-a_i+r_i(r_i+1)/2$.
On note 
$$L=Ri_{d-r_c,*}w_{>a_c}i_{d-r_c}^*\dots Ri_{d-r_1,*}w_{>a_1}i_{d-r_1}^*Rj_*\F^{\K}(V)$$
et $i=i_{g,r}$.
Soit $S=\{r_1,\dots,r_c,r\}$.

Alors on a un isomorphisme canonique
$$i^*L\simeq\bigoplus_C L_C,$$
où $C$ parcourt l'ensemble des doubles classes de $\Pa_S(\Q)\QP_r(\Af)\sous\Pa_r(\Q)\QP_r(\Af)/\Hr_r$,
et de plus, si $h$ est un représentant dans $\Pa_r(\Z)\QP_r(\widehat{\Z})$
de la double classe $C$, on a un isomorphisme
$$L_C\simeq\F^{\K_r}R\Gamma(\Gamma_{S},R\Gamma(Lie(\N_S),hg.V)_{<t_1,\dots,<t_{c}}).$$

\end{proposition}

\begin{proof} Si $r<r_c$, on pose $c'=c+1$, $r_{c'}=r$, $a_{c'}=-\infty$ et $t_{c'}=+\infty$
; si $r=r_c$, on pose $c'=c$.
On a donc $S=\{r_1,\dots,r_{c'}\}$,
avec $r_1>\dots>r_{c'}=r$.
D'après la définition de $L$, on a 
une décomposition 
$$i^*L=\bigoplus L_{X_1,\dots,X_{c'}},$$
où $(X_1,\dots,X_{c'})$ parcourt les $c'$-uplets de strates de $M^{\K}(\G,\X)^*$
de la forme $X_i=Im(i_{b_i,r_i})$
tels que $X_{i+1}\subset\overline{X}_i$ pour $1\leq i\leq c'-1$
et $X_{c'}=Im(i_{g,r})$,
et
$$L_{X_1,\dots,X_{c'}}=w_{>a_{c'}}i^*Ri_{b_{c'-1},r_{c'-1}*}w_{>a_{c'-1}}i_{b_{c'-1},r_{c'-1}}^*\dots Ri_{b_1,r_1*}w_{>a_1}i_{b_1,r_1}^*Rj_*\F^{\K}V.$$
Soit $(X_1,\dots,X_{c'})$ un tel $c'$-uplet.
On choisit $g_1\in\G(\widehat{\Z})$ tel que $X_1=Im(i_{g_1,r_1})$.
Pour tout $i\in\{1,\dots,c'-1\}$,
$\xymatrix@C=50pt{X_{i+1}\subset\overline{X}_i& M^{\K_{r_i}}(\G_{r_i},\X_{r_i})^*\ar[l]^-{\sim}_-{\overline{i}_{g_i\dots g_1,r_i}}}$
est la strate correspondant au sous-groupe parabolique maximal 
$\RP_{i+1}=(\Pa_{r_{i+1}}\cap\QP_{r_i})/\N_{r_i}$ de $\G_{r_i}$ et à un
$h_{i+1}\in\G_{r_i}(\widehat{\Z})$;
on choisit $g_{i+1}\in\QP_{r_{i}}(\widehat{\Z})$ relevant $h_{i+1}$.\newline
Comme $X_{c'}=Im(i_{g,r})$, il existe $h\in\Pa_r(\Q)\QP_r(\Af)$ tel
que $hg\K=g_{c'}\dots g_1\K$, 
et $h$ est dans $\Pa_r(\Z)\QP_r(\widehat{\Z})=\Pa_r(\Q)\QP_r(\Af)\cap\G(\widehat{\Z})$
car $g,g_1,\dots,g_{c'}\in\G(\widehat{\Z})$ et $\K\subset\G(\widehat{\Z})$.\newline
Notons, pour tout $i\in\{1,\dots,c'\}$, $S_i=\{r_1,\dots,r_i\}$.\newline
Soit $i\in\{1,\dots,c'-1\}$.
Le sous-groupe distingué $\QP_{R_{i+1}}$ de $\RP_{i+1}$ défini par Pink
est $\QP_{R_{i+1}}=\QP_{r_{i+1}}/(\QP_{r_{i+1}}\cap\N_{r_i})$,
le radical unipotent de $\RP_{i+1}$
est $\N_{R_{i+1}}=\N_{r_{i+1}}/(\N_{r_{i+1}}\cap\N_{r_i})=\N_{r_{i+1}}\N_{r_i}/\N_{r_i}$,
on a $\QP_{R_{i+1}}/\N_{R_{i+1}}=\G_{r_{i+1}}$,
et le quotient de Levi $\Le_{R_{i+1}}=\RP_{i+1}/\N_{R_{i+1}}$ 
s'écrit $\Le_{R_{i+1}}=\Le_{\lin,R_{i+1}}\times\G_{r_{i+1}}$,
avec $\Le_{\lin,R_{i+1}}=(\Pa_{r_{i+1}}\cap\QP_{r_i})/\QP_{r_{i+1}}\N_{r_i}$.
On voit facilement que
$\Le_{\ell,S_{i+1}}=\Le_{\lin,S_i}\times\Le_{\lin,R_{i+1}}$
et que $\N_{R_{i+1}}=\N_{S_{i+1}}/\N_{S_i}$.
On en déduit en particulier que le groupe
$$\begin{array}{rcl}\Gamma_{R_{i+1}} & = & (\K_{r_i}\cap\Le_{\lin,R_{i+1}}(\Q)\N_{R_{i+1}}(\Af))/(\K_{r_i}\cap\N_{R_{i+1}}(\Af)) \\
& = & Ker(\Le_{\lin,R_{i+1}}(\Z)\fl \Le_{\lin,R_{i+1}}(\Z/n\Z))\end{array}$$
vérifie $\Gamma_{S_{i+1}}=\Gamma_{S_i}\times\Gamma_{R_{i+1}}$.
Comme $\Le_{\lin,S_i}\subset\Le_{\lin,r_i}$ et $\G_{r_i}$ commutent
et que $\Se_1,\dots,\Se_{r_i}$ agissent trivialement sur $\N_{R_{i+1}}$,
on a un isomorphisme canonique
\begin{flushleft}$\displaystyle{R\Gamma(\Gamma_{R_{i+1}},R\Gamma(Lie(\N_{R_{i+1}}),g_{i+1}.R\Gamma(\Gamma_{S_i},R\Gamma(Lie(\N_{S_i}),g_i\dots g_1.V)_{<t_1,\dots,<t_{i}}))_{<t_{{i+1}}})\simeq}$\end{flushleft}
\begin{flushright}$\displaystyle{R\Gamma(\Gamma_{S_{i+1}},R\Gamma(Lie(\N_{S_{i+1}}),g_{i+1}\dots g_1.V)_{<t_1,\dots,<t_{{i+1}}}),}$\end{flushright}
d'où, grâce au théorème de Pink et au lemme \ref{relation_tronques}, 
un isomorphisme canonique
\begin{flushleft}$\displaystyle{w_{>a_{{i+1}}}i_{g_{i+1}\dots g_1,r_{i+1}}^*Ri_{g_i\dots g_1,r_i*}\F^{\K_{r_i}}R\Gamma(\Gamma_{S_i},R\Gamma(Lie(\N_{S_i}),g_i\dots g_1.V)_{<t_1,\dots,<t_{i}})\simeq}$\end{flushleft}
\begin{flushright}$\displaystyle{\F^{\K_{r_{i+1}}}R\Gamma(\Gamma_{S_{i+1}},R\Gamma(Lie(\N_{S_{i+1}}),g_{i+1}\dots g_1.V)_{<t_1,\dots,<t_{{i+1}}}).}$\end{flushright}
Une récurrence sur $i$ donne alors un isomorphisme canonique
$$\begin{array}{rcl}L_{X_1,\dots,X_{c'}} & \simeq & \displaystyle{\F^{\K_r}R\Gamma(\Gamma_S,R\Gamma(Lie(\N_S),g_{c'}\dots g_1.V)_{<t_1,\dots,<t_{c'}})} \\
& \simeq & \displaystyle{\F^{\K_r}R\Gamma(\Gamma_S,R\Gamma(Lie(\N_S),hg.V)_{<t_1,\dots,<t_{c'}})} \\
& = & \displaystyle{\F^{\K_r}R\Gamma(\Gamma_S,R\Gamma(Lie(\N_S),hg.V)_{<t_1,\dots,<t_c})}.\end{array}$$
Il reste à compter les $c'$-uplets $(X_1,\dots,X_{c'})$.
On a associé à $(X_1,\dots,X_{c'})$ (de manière non unique) un $(c'+1)$-uplet $(g_1,\dots,g_{c'},h)$,
avec $g_1\in\G(\widehat{\Z})$, $g_i\in\QP_{r_{i-1}}(\widehat{\Z})$ pour $i\geq 2$ et $h\in\Pa_r(\Z)\QP_r(\widehat{\Z})$
tels que $hg\K=g_{c'}\dots g_1\K$.
On note $(X_1,\dots,X_{c'})=\underline{X}_{g_1,\dots,g_{c'},h}$.\newline
On a $\underline{X}_{g_1,\dots,g_{c'},h}=\underline{X}_{g'_1,\dots,g'_{c'},h'}$ si et seulement si
$$\Pa_{r_1}(\Q)\QP_{r_1}(\Af)g_1\K=\Pa_{r_1}(\Q)\QP_{r_1}(\Af)g'_1\K$$
$$\Pa_{\{r_1,r_2\}}(\Q)\QP_{r_2}(\Af)g_2g_1\K=\Pa_{\{r_1,r_2\}}(\Q)\QP_{r_2}(\Af)g'_2g'_1\K$$
$$\dots$$
$$\Pa_{\{r_1,\dots,r_{c'}\}}(\Q)\QP_{r_{c'}}(\Af)g_{c'}\dots g_1\K=\Pa_{\{r_1,\dots,r_{c'}\}}(\Q)\QP_{r_{c'}}(\Af)g'_{c'}\dots,g'_1\K.$$
Il est évident que chaque ligne implique la précédente, 
et que la dernière condition est équivalente à
$$\Pa_S(\Q)\QP_r(\Af)hg\K=\Pa_S(\Q)\QP_r(\Af)h'g\K,$$
elle-même équivalente à
$$\Pa_S(\Q)\QP_r(\Af)h\Hr_r=\Pa_S(\Q)\QP_r(\Af)h'\Hr_r,$$
car $g\K g^{-1}=\K$ ($g\in\G(\widehat{\Z})$).\newline
On en déduit finalement que l'application
$\underline{X}_{g_1,\dots,g_{c'},h}\fle\Pa_S(\Q)\QP_r(\Af)h\Hr_r$
induit une bijection entre l'ensemble des $c'$-uplets $(X_1,\dots,X_{c'})$
et le double quotient\newline
$\Pa_S(\Q)\QP_r(\Af)\sous\Pa_r(\Q)\QP_r(\Af)/\Hr_r$.

\end{proof}

\begin{remarque} Nous pouvons maintenant rendre plus explicite le rapport entre les complexes pondérés
définis ici
et ceux de [GHM].
Soient $t_0,\dots,t_{d-1}\in\Z\cup\{\pm\infty\}$
et $V$ une représentation algébrique de $\G$,
qu'on suppose pure de poids $0$ pour simplifier.
Pour tout $r\in\{1,\dots,d\}$, on note
 $a_r=-t_{d-r}+(d-r)(d+1-r)/2$.
Alors, d'après le théorème \ref{th:simplification_formule_traces},
on a une égalité dans le groupe de Grothendieck de $D^b_m(M^{\K}(\G,\X)^*,\Q_\ell)$
$$[W^{\geq t_0,\dots,\geq t_{d-1}}V(\Q_\ell)]=\sum_{1\leq r_1<\dots<r_c\leq d}(-1)^c[Ri_{r_c*}w_{>a_{r_c}}i_{r_c}^*\dots Ri_{r_1*}w_{>a_1}i_{r_1}^*Rj_*\F^{\K}V(\Q_\ell)].$$
Or, d'après le calcul explicite des
$$L_{r_1,\dots,r_c}=Ri_{r_c*}w_{>a_{r_c}}i_{r_c}^*\dots Ri_{r_1*}w_{>a_{r_1}}i_{r_1}^*Rj_*\F^{\K}V(\Q_\ell)$$
qui a été fait dans la proposition ci-dessus,
il existe une manière naturelle de relever ces complexes en des complexes
sur $\Mod_{d,n}^*[1/\ell]$, qu'on notera encore $L_{r_1,\dots,r_c}$.
Notons $L_{r_1,\dots,r_c}(\C)$ le complexe de faisceaux de $\Q_\ell$-espaces vectoriels
sur $\Mod_{d,n}^*(\C)$ déduit du complexe
$L_{r_1,\dots,r_c}$ sur $\Mod_{d,n}^*[1/\ell]$.
Alors la somme alternée
$$\sum_{1\leq r_1<\dots<r_c\leq d}(-1)^c[L_{r_1,\dots,r_c}(\C)]$$
est égale à la classe (dans le groupe de Grothendieck
de la catégorie dérivée de la catégorie des faisceaux de $\Q_\ell$-espaces vectoriels
sur $\Mod_{d,n}^*(\C)$)
de l'image directe par le morphisme canonique de la compactification de Borel-Serre réductive
de $M^{\K}(\G,\X)(\C)$ sur $\Mod_{d,n}^*(\C)=M^{\K}(\G,\X)^*(\C)$
du complexe pondéré de [GHM] associé au profil de poids $(t_0,\dots,t_{d-1})$
et à coefficients dans $V$.

\end{remarque}

\section{Correspondances de Hecke}

Dans toute cette partie, les flèches marquées ``CB'' seront des flèches
de changement de base.

\subsection{Correspondances cohomologiques et troncature par le poids}

Soient $\overline{c}_1:X'\fl X_1,\overline{c}_2:X'\fl X_2$ deux morphismes finis
entre des schémas séparés de type fini sur $\Fi_q$.

\begin{definition} Soient $j_1:Y_1\fl X_1,j_2:Y_2\fl X_2,j':Y'\fl X'$ des immersions 
localement fermées.
On suppose que $\overline{c}_1(Y')\subset Y_1$ et $\overline{c}_2(Y')\subset Y_2$,
et on note $c_1:Y'\fl Y_1$ et $c_2:Y'\fl Y_2$ les morphismes obtenus par restriction.

Alors, pour tous $K\in D^b_c(Y_1,\Q_\ell),L\in D^b_c(Y_2,\Q_\ell)$ et
pour toute correspondance cohomologique $u:c_1^*K\fl c_2^!L$ de $K$ à $L$
à support dans $(c_1,c_2)$, on appelle image de $u$ par $(j_1,j_2)$
la correspondance cohomologique de $Rj_{1*}K$ à $Rj_{2*}L$ à support dans 
$(\overline{c}_1,\overline{c}_2)$ suivante :
$$\overline{c}_1^*Rj_{1*}K\stackrel{CB}{\fl} Rj'_*c_1^*K\stackrel{u}{\fl} Rj'_*c_2^!L\stackrel{CB}{\fl}\overline{c}_2^!Rj_{2*}L.$$

Si $X_1=X_2$ et $j_1=j_2=j$, on note aussi $\overline{u}=Rj_*u$.

\end{definition}

\begin{lemme}\label{lemme:unicite_prolong}(\cite{F} 1.3.1) Supposons que $X_1=X_2$ et que $j=j_1=j_2$ et $j'$ sont des immersions ouvertes.
Alors, pour toute $u:c_1^*K\fl c_2^!L$,
 $Rj_*u$ est l'unique prolongement de $u$ en une correspondance cohomologique
$\overline{c}_1^*Rj_*K\fl\overline{c}_2^!Rj_*L$.

\end{lemme}

Dans la suite, on suppose que $X_1=X_2=X$.
Soient $(S_k)_{0\leq k\leq n}$ et $(S'_k)_{0\leq k\leq n}$ des stratifications de $X$ et $X'$
(au sens de la définition \ref{definition_stratification}).
On note $i_k:S_k\fl X$ (resp. $i'_k:S'_k\fl X'$) l'inclusion, $U=S_0$ (resp. $U'=S'_0$) 
et $j=i_0$ (resp. $j'=i'_0$).

On suppose que, pour tout $k\in\{0,\dots,n\}$, $\overline{c}_1$ et $\overline{c}_2$
envoient $S'_k$ dans $S_k$. 

On note $c_1^k,c_2^k:S'_k\fl S_k$ les morphismes obtenus;
si $k=0$, on note aussi $c_1=c_1^0$ et $c_2=c_2^0$.

\begin{lemme}\label{lemme:corr_ponderees} Soient $K,L\in D^b_m(X,\Q_\ell)$, $u:\overline{c}_1^*K\fl\overline{c}_2^!L$
une correspondance cohomologique à support dans $(\overline{c}_1,\overline{c}_2)$
et $\as\in (\Z\cup\{\pm\infty\})^{n+1}$.
Alors il existe une et une seule correspondance cohomologique
$w_{\leq\as}u:\overline{c}_1^*w_{\leq\as}K\fl \overline{c}_2^!w_{\leq\as}L$
(resp. $w_{>\as}u:\overline{c}_1^*w_{>\as}K\fl\overline{c}_2^!w_{>\as}L$)
qui fait commuter le diagramme
$$\xymatrix{\overline{c}_1^*w_{\leq\as}K\ar[r]\ar[d]^{w_{\leq\as}u} & \overline{c}_1^*K\ar[d]^{u} \\
\overline{c}_2^!w_{\leq\as}L\ar[r] & \overline{c}_2^!L}$$
$$\left(resp.\qquad\vcenter{\xymatrix{\overline{c}_1^*K\ar[r]\ar[d]^{u} & \overline{c}_1^*w_{>\as}K\ar[d]^{w_{>\as}u} \\
\overline{c}_2^!L\ar[r] & \overline{c}_2^!w_{>\as}L}}\right)$$

\end{lemme}

\begin{proof} Il suffit de traiter le cas de $w_{\leq\as}$, car celui de $w_{>\as}$
s'en déduit par dualité.
Considérons le diagramme suivant, dont les lignes sont distinguées
$$\xymatrix{\overline{c}_1^*w_{\leq\as}K\ar[r] & \overline{c}_1^*K\ar[r]\ar[d]_u & \overline{c}_1^*w_{>\as}K\ar[r]^{+1} & \\
\overline{c}_2^!w_{\leq\as}L\ar[r] & \overline{c}_2^!L\ar[r] & \overline{c}_2^!w_{>\as}L\ar[r]^{+1} &}$$
D'après le (iv) de la proposition \ref{prop2:t_structures_recollees},
$\overline{c}_1^*w_{\leq\as}K\in\DP^{\leq\as}(X')$ et $\overline{c}_2^!w_{>\as}L\in\DP^{>\as}(X')$.
En appliquant le (i) de la même proposition, on trouve
$R\Hom(\overline{c}_1^*w_{\leq\as}K,\overline{c}_2^!w_{>\as}L)=0$,
donc le morphisme canonique
$$R\Hom(\overline{c}_1^*w_{\leq\as}K,\overline{c}_2^!w_{\leq\as}L)\fl R\Hom(\overline{c}_1^*w_{\leq\as}K,\overline{c}_2^!L)$$
est un isomorphisme, 
ce qui donne l'existence et l'unicité de $w_{\leq\as}u$.

\end{proof}

Nous allons expliciter les correspondances $w_{\leq\as}u$ et $w_{>\as}u$ dans deux cas particuliers.
On aura besoin d'un morphisme fonctoriel déduit de l'isomorphisme de changement de base propre.
Considérons un diagramme cartésien aux éléments nilpotents près
$$\xymatrix{Y'\ar[r]^{f'}\ar[d]_{g'} & X'\ar[d]^{g} \\
Y\ar[r]_{f} & X}$$
On a un isomorphisme de changement de base $Rf'_*{g'}^!\iso g^!Rf_*$, d'où un morphisme fonctoriel
$${f'}^*g^!\stackrel{adj}{\fl}{f'}^*g^!Rf_*f^*\stackrel{\sim}{\longleftarrow}{f'}^*Rf'_*{g'}^!f^*\stackrel{adj}{\fl}{g'}^!f^*.$$
Si $j:U\fl Y'$ est un morphisme étale, on en déduit un morphisme fonctoriel
$$(f'j)^*g^!=j^!{f'}^*g^!\fl j^!{g'}^!f^*=(g'j)^!f^*.$$

Revenons au calcul de $w_{\leq\as}u$ et $w_{>\as}u$.

\begin{lemme}\label{lemme:restr_corr_bord} Soient $a\in\Z\cup\{\pm\infty\}$ et $k\in\{0,\dots,n\}$.
On définit $\as,\as'\in (\Z\cup\{\pm\infty\})^{n+1}$ par $a_r=-a'_r=+\infty$ si $r\not=k$
et $a_k=a'_k=a$.
On suppose que $c_1^k$ et $c_2^k$ sont étales.
Alors, pour tous $K,L\in D^b_m(X,\Q_\ell)$ et
pour toute correspondance cohomologique $u:\overline{c}_1^*K\fl\overline{c}_2^!L$,
$w_{>\as}u$ est égal au composé
$$\xymatrix@C=15pt@R=10pt{\overline{c}_1^*w_{>\as}K\ar@{=}[r] & \overline{c}_1^*Ri_{k*}w_{>a}i_k^*K\ar[r]^{CB} & Ri'_{k*}c_1^{k*}w_{>a}i_k^*K\ar@{=}[r] & Ri'_{k*}w_{>a}c_1^{k*}i_k^*K\ar@{=}[r] & Ri'_{k*}w_{>a}{i'_k}^*\overline{c}_1^*K\ar[d]^{u} \\
\overline{c}_2^!w_{>\as}L\ar@{=}[r] & \overline{c}_2^!Ri_{k*}w_{>a}i_k^*L & Ri'_{k*}c_2^{k!}w_{>a}i_k^*L\ar[l]_{CB} & Ri'_{k*}w_{>a}c_2^{k!}i_k^*L\ar@{=}[l] & R'_{k*}w_{>a}{i'_k}^*\overline{c}_2^!L\ar[l]}$$
et $w_{\leq\as'}u$ est égal au composé
$$\xymatrix@C=15pt@R=10pt{\overline{c}_1^*w_{\leq\as'}K\ar@{=}[r] & \overline{c}_1^*i_{k!}w_{\leq a}i_k^!K\ar[r]^{CB} & i'_{k!}c_1^{k*}w_{\leq a}i_k^!K\ar@{=}[r] & i'_{k!}w_{\leq a}c_1^{k*}i_k^!K\ar[r] & i'_{k!}w_{\leq a}{i'_k}^!\overline{c}_1^*K\ar[d]^{u} \\
\overline{c}_2^!w_{\leq\as'}L & \overline{c}_2^!i_{k!}w_{\leq a}i_k^!L\ar@{=}[l] & i'_{k!}c_2^{k!}w_{\leq a}i_k^!L\ar[l]_{CB} & i'_{k!}w_{\leq a}c_2^{k!}i_k^!L\ar@{=}[l] & i'_{k!}w_{\leq a}{i'_k}^!\overline{c}_2^!L\ar@{=}[l]}$$

\end{lemme}

\begin{proof} Il suffit de traiter le cas de $w_{>\as}u$,
l'autre cas s'en déduisant par dualité.
Notons $v:\overline{c}_1^*w_{>\as}K\fl \overline{c}_2^!w_{>\as}L$ le morphisme défini dans l'énoncé.
D'après le lemme \ref{lemme:corr_ponderees}, il suffit de montrer que le diagramme suivant est commutatif
$$\xymatrix{\overline{c}_1^*K\ar[d]_{u}\ar[r] & \overline{c}_1^*w_{>\as}K\ar[d]^{v} \\
\overline{c}_2^!L\ar[r] & \overline{c}_2^!w_{>\as}L}$$
Le morphisme $v$ est le composé de trois morphismes :
\begin{itemize}
\item[(1)] un morphisme $v_1:\overline{c}_1^*w_{>\as}K\fl Ri'_{k*}w_{>a}{i'_k}^*\overline{c}_1^*K=w_{>\as}\overline{c}_1^*K$;
\item[(2)] un morphisme $w_{>\as}\overline{c}_1^*K\fl w_{>\as}\overline{c}_2^!L$,
obtenu en appliquant à $u$ le foncteur $w_{>\as}$;
\item[(3)] un morphisme $v_2:w_{>\as}\overline{c}_2^!L=Ri'_{k*}w_{>a}{i'_k}^*\overline{c}_2^!L\fl \overline{c}_2^!w_{>\as}L$.

\end{itemize}
Il suffit de montrer que les deux diagrammes suivants commutent
$$\xymatrix{\overline{c}_1^*K\ar[r]\ar[rd] & \overline{c}_1^*w_{>\as}K\ar[d]^{v_1} & \overline{c}_2^!L\ar[r]\ar[rd] & w_{>\as}\overline{c}_2^!L\ar[d]^{v_2} \\
& w_{>\as}\overline{c}_1^*K & & \overline{c}_2^!w_{>\as}L}$$
Le premier diagramme se réécrit
$$\xymatrix{\overline{c}_1^*K\ar[r]^-{adj}\ar[rdd]_-{CB}\ar@{}[rd]|{(1)} & \overline{c}_1^*Ri_{k*}i_k^*K\ar[d]^{CB}\ar[r]\ar@{}[rd]|{(2)} & \overline{c}_1^*Ri_{k*}w_{>a}i_k^*K\ar[d]^{CB} \\
& Ri'_{k*}c_1^{k*}i_k^*K\ar[r]\ar@{=}[d]\ar@{}[rd]|{(3)} & Ri'_{k*}c_1^{k*}w_{>a}i_k^*K\ar[d]^{\wr} \\
 & Ri'_{k*}{i'_k}^*\overline{c}_1^*K\ar[r] & Ri'_{k*}w_{>a}{i'_k}^*\overline{c}_1^*K }$$
$(1)$ est commutatif par des arguments standard,
$(2)$ et $(3)$ sont clairement commutatifs.
Pour prouver que le deuxième diagramme est commutatif, 
on se ramène de la même manière à montrer la commutativité du diagramme
$$\xymatrix{\overline{c}_2^!L\ar[r]^-{adj}\ar[d]_{adj} & Ri'_{k*}{i'_k}^*\overline{c}_2^!L\ar[d] \\
\overline{c}_2^!Ri_{k*}i_k^*L & Ri'_{k*}c_2^{k!}i_k^*L\ar[l]_{CB}}$$
Si on remplace la flèche de droite par sa définition, ce diagramme se réécrit
$$\xymatrix{\overline{c}_2^!L\ar[r]^-{adj}\ar[d]_{adj} & Ri'_{k*}{i'_k}^*\overline{c}_2^!L\ar[r]^-{adj} & Ri'_{k*}{i'_k}^*\overline{c}_2^!L \\
\overline{c}_2^!Ri_{k*}i_k^*L & Ri'_{k*}{c_2^k}^!i_k^*L\ar[l]_-{CB} & Ri'_{k*}{i'_k}^*Ri'_{k*}{c_2^k}^!i_k^*L\ar[l]_-{adj}\ar[u]_{CB}^\wr}$$
En développant les flèches de changement de base, on obtient le diagramme ci-dessous
(où toutes les flèches sont des flèches d'adjonction), 
qui est clairement commutatif
$$\xymatrix@C=20pt{\overline{c}_2^!L\ar[r]\ar[d] & Ri'_{k*}{i'_k}^*\overline{c}_2^!L\ar[r] & Ri'_{k*}{i'_k}^*\overline{c}_2^!Ri_{k*}i_k^*L & Ri'_{k*}{i'_k}^*\overline{c}_2^!Ri_{k*}{c_2^k}_*{c_2^k}^!i_k^*L\ar[l]\ar@{=}[d] \\
\overline{c}_2^!Ri_{k*}i_k^*L\ar[urr] & & & Ri'_{k*}{i'_k}^*\overline{c}_2^!\overline{c}_{2*}Ri'_{k*}{c_2^k}^!i_k^*L \\
\overline{c}_2^!Ri_{k*}{c_2^k}_*{c_2^k}^!i_k^*L\ar[u]\ar@{=}[r]\ar[rrruu] & \overline{c}_2^!\overline{c}_{2*}Ri'_{k*}{c_2^k}^!i_k^*L\ar[rru] & Ri'_{k*}{c_2^k}^!i_k^*L\ar[r]\ar[l] & Ri'_{k*}{i'_k}^*Ri'_{k*}{c_2^k}^!i_k^*L\ar[u]\ar@/_1.5pc/[l] }$$

\end{proof}

Nous appelons ``groupe de Grothendieck des correspondances cohomologiques à support
dans $(\overline{c}_1,\overline{c}_2)$'' le groupe engendré par les classes
d'isomorphisme de triplets $(K,L,u)$, où $K,L\in D^b_m(X,\Q_\ell)$ et $u$ est un morphisme
$\overline{c}_1^*K\fl \overline{c}_2^!L$
(un isomorphisme $(K,L,u)\iso (K',L',u')$ est un couple d'isomorphismes $(f:K\iso K',g:L\iso L')$ 
tel que $\overline{c}_2^!(g)\circ u=u'\circ\overline{c}_1^*(f)$),
soumis aux relations : 
$[(K',L',u')]=[(K,L,u)]+[(K'',L'',u'')]$
s'il existe des triangles distingués $K\fl K'\fl K''\stackrel{+1}{\fl}$
et $L\fl L'\fl L''\stackrel{+1}{\fl}$ tels que le diagramme suivant commute
$$\xymatrix{\overline{c}_1^*K\ar[r]\ar[d]_{u} & \overline{c}_1^*K'\ar[r]\ar[d]_{u'} & \overline{c}_1^*K''\ar[r]^{+1}\ar[d]_{u''} & \\
\overline{c}_2^!L\ar[r] & \overline{c}_2^!L'\ar[r] & \overline{c}_2^!L''\ar[r]^{+1} & }$$

La proposition suivante se prouve exactement comme le théorème \ref{th:simplification_formule_traces}.

\begin{proposition}\label{prop:s_f_t} Pour tout $\as\in (\Z\cup\{\pm\infty\})^{n+1}$, 
pour tous $K,L\in \DP^{\leq a_0}(U)$ et
pour toute correspondance cohomologique
 $u:c_1^*K\fl c_2^!L$, 
on a
$$[w_{\leq\as}Rj_*u]=
\sum_{1\leq n_1<\dots n_r\leq n}(-1)^r[Ri_{n_r*}w_{>a_r}i_{n_r}^*\dots Ri_{n_1*}w_{>a_1}i_{n_1}^*Rj_*u]$$
dans le groupe de Grothendieck des correspondances cohomologiques à support dans $(\overline{c}_1,\overline{c}_2)$.

\end{proposition}

\subsection{Correspondances de Hecke sur les complexes pondérés}

Dans cette section, $(\G,\X)$ est la donnée de Shimura $(\GSp_{2d,\Q},\X_d)$
de la section 1.1,
et on utilise les notations de la partie 4.
On fixe deux entiers $n,m\geq 3$ tels que $n$ divise $m$, et un nombre premier $p\not=\ell$ 
qui ne divise pas $m$. On note $\K=\K_d(n)$ et $\K'=\K_d(m)$ et,
comme dans la partie 4, on travaille sur les réductions modulo $p$
des variétés de Shimura.

\begin{definition}
Soient $g\in\G(\Af^p)$ tel que $g^{-1}\K'g\subset\K$ et
$V\in D^b(Rep_{\G(\Q_\ell)})$. 
La multiplication par $g$ induit un isomorphisme $g.V\iso V$,
d'où une correspondance cohomologique sur $\F^{\K}V$ à support dans $(T_g,T_1)$,
$$u_g:T_g^*\F^{\K}V\simeq\F^{\K'}g.V\iso\F^{\K'}V\simeq T_1^*\F^{\K}V,$$
qu'on appelle correspondance de Hecke associée à $g$.

\end{definition}

Dans cette section, nous allons calculer, pour $V\in D^b(Rep_{\G(\Q_\ell)})$,
$0\leq r_c<\dots<r_1\leq d-1$
et $a_1,\dots,a_c\in\Z\cup\{\pm\infty\}$, la correspondance
$Ri_{d-r_c,*}w_{>a_c}i_{d-r_c}^*\dots Ri_{d-r_1,*}w_{>a_1}i_{d-r_1}^*Rj_*u_g$
(d'après la proposition \ref{prop:s_f_t}, on pourra en déduire la classe dans le groupe de Grothendieck
d'une correspondance $w_{\leq\as}Rj_*u_g$).
Pour tout $s\in\{1,\dots,c\}$, on pose $t_s=-a_s+r_s(r_s+1)/2$.

Notons $S=\{r_1,\dots,r_c\}$ et $\mathcal{E}$ l'ensemble des $c$-uplets $(X_1,\dots,X_c)$
tels que $X_i\subset M^{\K}(\G,\X)^*$ soit de la forme $Im(i_{b_i,r_i})$, $b_i\in\G(\widehat{\Z})$,
et $X_{i+1}\subset\overline{X}_i$ pour $1\leq i\leq c-1$.

On définit une application
$$\left\{\begin{array}{rcl}\G(\widehat{\Z})\times\QP_{r_1}(\widehat{\Z})\times\dots\times\QP_{r_{c-1}}(\widehat{\Z}) & \fl & \mathcal{E} \\
(g_1,\dots,g_c) & \fle & \underline{X}_{g_1,\dots,g_c}\end{array}\right.$$
de la manière suivante : $\underline{X}_{g_1,\dots,g_c}=(X_1,\dots,X_c)$
avec $X_1=Im(i_{g_1,r_1})$, et, pour tout $\protect{j\in\{1,\dots,c-1\}}$,
$X_{j+1}$ la strate de bord de 
$\xymatrix@C=40pt{M^{\K_{r_j}}(\G_{r_j},\X_{r_j})^*\ar[r]^-{\overline{i}_{g_j\dots g_1,r_j}}_-{\sim} & \overline{X}_j}$
associée au sous-groupe parabolique maximal
$(\Pa_{r_{j+1}}\cap\QP_{r_j})/\N_{r_j}$ de $\G_{r_j}$
et à l'image de $g_{j+1}$ dans $\G_{r_j}(\widehat{\Z})$.

Cette application est surjective, et on a $\underline{X}_{g_1,\dots,g_c}=\underline{X}_{g'_1,\dots,g'_c}$ si et seulement si
$$\Pa_{r_1}(\Q)\QP_{r_1}(\Af)g_1\K=\Pa_{r_1}(\Q)\QP_{r_1}(\Af)g'_1\K$$
$$\Pa_{\{r_1,r_2\}}(\Q)\QP_{r_2}(\Af)g_2g_1\K=\Pa_{\{r_1,r_2\}}(\Q)\QP_{r_2}(\Af)g'_2g'_1\K$$
$$\dots$$
$$\Pa_{\{r_1,\dots,r_c\}}(\Q)\QP_{r_c}(\Af)g_c\dots g_1\K=\Pa_{\{r_1,\dots,r_c\}}(\Q)\QP_{r_c}(\Af)g'_c\dots g'_1\K.$$
Comme chaque ligne implique la précédente, on en déduit que $\underline{X}_{g_1,\dots,g_c}=\underline{X}_{g_c\dots g_1,1,\dots,1}$,
et que $g\fle\underline{X}_{(g,1,\dots,1)}$
induit une bijection $\varphi:\Pa_S(\Q)\QP_{r_c}(\Af)\sous\G(\Af)/\K\iso\mathcal{E}$.

D'après la proposition \ref{restr_tronques_bord}, on a un isomorphisme canonique
$$Ri_{d-r_c,*}w_{>a_c}i_{d-r_c}^*\dots Ri_{d-r_1,*}w_{>a_1}i_{d-r_1}^*Rj_*\F^{\K}V
\simeq\bigoplus_{C\in\Pa_S(\Q)\QP_{r_c}(\Af)\sous\G(\Af)/\K}L_{C},$$
où, pour tout $C$, si $\varphi(C)=(Im(i_{b_1,r_1}),\dots,Im(i_{b_c,r_c}))\in\mathcal{E}$, alors
$$L_{C}=Ri_{b_c,r_c,*}w_{>a_c}i_{b_c,r_c}^*\dots Ri_{b_1,r_1,*}w_{>a_1}i_{b_1,r_1}^*Rj_*\F^{\K}V.$$

La correspondance $Ri_{d-r_c,*}w_{>a_c}i_{d-r_c}^*\dots Ri_{d-r_1,*}w_{>a_1}i_{d-r_1}^*Rj_*u_g$
est donc donnée par une matrice 
$(u_{C_1,C_2})_{C_1,C_2\in\Pa_S(\Q)\QP_{r_c}(\Af)\sous\G(\Af)/\K}$,
dont nous allons calculer les coefficients.

Dans la suite, on utilisera les notations de la section 1.2 
(en particulier $\Hr_r$, $\Hr_{\ell,r}$, etc)
et on notera avec un $'$ les groupes analogues obtenus en remplaçant $\K$ par $\K'$.

Soient $C_1,C_2\in\Pa_S(\Q)\QP_{r_c}(\Af)\sous\G(\Af)/\K$.
On fixe $h_1\in C_1\cap\G(\widehat{\Z})$ et $h_2\in C_2\cap\G(\widehat{\Z})$.
D'après la proposition \ref{restr_tronques_bord}, on a des isomorphismes
$$L_{C_1}\simeq Ri_{h_1,r_c,*}\F^{\K_{r_c}}R\Gamma(\Gamma_S,R\Gamma(Lie(\N_S),h_1.V)_{<t_1,\dots,<t_c})$$
$$L_{C_2}\simeq Ri_{h_2,r_c,*}\F^{\K_{r_c}}R\Gamma(\Gamma_S,R\Gamma(Lie(\N_S),h_2.V)_{<t_1,\dots,<t_c}).$$

Soient $C'$ une double classe dans $\Pa_S(\Q)\QP_{r_c}(\Af)\sous\G(\Af)/\K'$
et $h\in C'\cap\G(\widehat{\Z})$.
On suppose que $C_1=C'g\K$ et $C_2=C'\K$.
Il existe donc $q_1,q_2\in\Pa_S(\Q)\QP_{r_c}(\Af)\cap\G(\Af^p)$ tels que
$q_1h_1\in hg\K$ et $q_2h_2\in h\K$, et les diagrammes suivants sont commutatifs
$$\xymatrix{M^{\K'_{r_c}}(\G_{r_c},\X_{r_c})\ar[r]^{i'_{h,{r_c}}}\ar[d]_{T_{\overline{q}_1}} & M^{\K'}(\G,\X)^*\ar[d]_{\overline{T}_g} & M^{\K'_{r_c}}(\G_{r_c},\X_{r_c})\ar[r]^{i'_{h,{r_c}}}\ar[d]_{T_{\overline{q}_2}} & M^{\K'}(\G,\X)^*\ar[d]_{\overline{T}_1} \\
M^{\K_{r_c}}(\G_{r_c},\X_{r_c})\ar[r]^{i_{h_1,{r_c}}} & M^{\K}(\G,\X)^* & M^{\K_{r_c}}(\G_{r_c},\X_{r_c})\ar[r]^{i_{h_2,{r_c}}} & M^{\K}(\G,\X)^*}$$
où $\overline{q}_1$ (resp. $\overline{q}_2$) est l'image de $q_1$ (resp. $q_2$)
dans $\G_{r_c}(\Af^p)$.

On note $v_{C'}$ la correspondance cohomologique de $L_{C_1}$ à $L_{C_2}$
à support dans $(\overline{T}_g,\overline{T}_1)$ image par $(i_{h_1,r_c},i_{h_2,r_c})$
de la correspondance
\begin{flushleft}$\displaystyle{T_{\overline{q}_1}^*i_{h_1,r_c}^*L_{C_1}\simeq\F^{\K'_{r_c}}
R\Gamma(\Gamma_S,R\Gamma(Lie(\N_S),q_1h_1.V)_{<t_1,\dots,<t_c})}$\end{flushleft}
\begin{flushright}$\displaystyle{\stackrel{q_2^{-1}h_2^{-1}q_1h_1}
{\overrightarrow{\hspace{2cm}}}\F^{\K'_{r_c}}R\Gamma(\Gamma_S,R\Gamma(Lie(\N_S),q_2h_2.V)
_{<t_1,\dots,<t_c})\simeq T_{\overline{q}_2}^*i_{h_2,r_c}^*L_{C_2},}$\end{flushright}
où la flèche du milieu est induite par l'isomorphisme
$q_1h_1.V\iso q_2h_2.V,v\fle q_{2,\ell}^{-1}h_{2,\ell}^{-1}q_{1,\ell}h_{1,\ell}.v$.

\begin{theoreme}\label{th:restr_corr_sigma} On a
$$u_{C_1,C_2}=[\Hr_{\ell,S}:\Hr'_{\ell,S}]\sum_{C'}v_{C'},$$
où la somme est sur les doubles classes $C'$ dans 
$\Pa_S(\Q)\QP_r(\Af)\sous\G(\Af)/\K'$
telles que $C_1=C'g\K$ et $C_2=C'\K$.

\end{theoreme}

\begin{proof} On note $\mathcal{E}'$ l'ensemble des $c$-uplets $(X'_1,\dots,X'_c)$
tels que, pour tout $i\in\{1,\dots,c\}$, $X'_i\subset M^{\K'}(\G,\X)^*$
soit de la forme $Im(i'_{b_i,r_i})$, $b_i\in\G(\widehat{\Z})$,
et, pour tout $i\in\{1,\dots,c-1\}$, $X'_{i+1}\subset \overline{X'}_i$.
On construit comme plus haut une bijection 
$\Pa_S(\Q)\QP_{r_c}(\Af)\sous\G(\Af)/\K'\iso\mathcal{E}'$
et, pour tout $C'\in\Pa_S(\Q)\QP_{r_c}(\Af)\sous\G(\Af)/\K'$, un complexe
$L'_{C'}$, qui, si $h\in C'\cap\G(\widehat{\Z})$,  est isomorphe à
$$\F^{\K'_{r_c}}R\Gamma(\Gamma_S,R\Gamma(Lie(\N_S),h.V)_{<t_1,\dots,<t_c}).$$

Les morphismes $\overline{T}_g$ et $\overline{T}_1$ induisent de manière évidente 
des applications $\mathcal{E}'\fl\mathcal{E}$, qu'on notera encore $\overline{T}_g$
et $\overline{T}_1$.
Si on identifie $\mathcal{E}$ (resp. $\mathcal{E}'$) à
$\Pa_S(\Q)\QP_{r_c}(\Af)\sous\G(\Af)/\K$ (resp. $\Pa_S(\Q)\QP_{r_c}(\Af)\sous\G(\Af)/\K'$),
$\overline{T}_g$ et $\overline{T}_1$ sont simplement les applications
$C'\fle C'g\K$ et $C'\fle C'\K$.

D'après le lemme \ref{lemme:restr_corr_bord}, $u_{C_1,C_2}$ s'écrit $\sum u_{C'}$,
pour $C'$ parcourant $\overline{T}_g^{-1}(C_1)\cap\overline{T}_1^{-1}(C_2)$,
où $u_{C'}$ est de la forme $\overline{T}_g^*L_{C_1}\fl L'_{C'}\fl\overline{T}_1^*L_{C_2}$.
Pour voir que $u_{C'}=[\Hr_{\ell,S}:\Hr'_{\ell,S}]v_{C'}$, 
il suffit de raisonner par récurrence sur $c$ comme dans la preuve de la
proposition \ref{restr_tronques_bord}, en utilisant la description du lemme
\ref{lemme:restr_corr_bord} et le corollaire de la proposition ci-dessous.

\end{proof}

\begin{proposition}\label{prop:suite_Pink} Soient $h,h'\in\G(\widehat{\Z})$,
$g\in\G(\Af^p)$ tel que $g^{-1}\K' g\subset\K$
et $r\in\{0,\dots,d-1\}$.
On note
$i=i_{h,r}:M_1=M^{\K_r}(\G_r,\X_r)\fl M^{\K}(\G,\X)^*$ 
(resp. $i'=i'_{h',r}:M_1'=M^{\K'_r}(\G_r,\X_r)\fl M^{\K'}(\G,\X)^*$).
On suppose que $\overline{T}_g(i'(M_1'))=i(M_1)$, c'est-à-dire que
$\Pa_r(\Q)\QP_r(\Af)h\K=\Pa_r(\Q)\QP_r(\Af)h'g\K$.
On a un diagramme commutatif
$$\xymatrix{M'_1\ar[d]_{d}\ar[r]^-{i'} & M^{\K'}(\G,\X)^*\ar[d]_{\overline{T}_g} & M^{\K'}(\G,\X)\ar[l]_-{j'}\ar[d]_{T_g} \\
M_1\ar[r]^-{i} & M^{\K}(\G,\X)^* & M^{\K}(\G,\X)\ar[l]_-{j}}$$
Le carré de droite est cartésien, et $T_g$ et $d$ sont finis étales.
Alors, pour tout $L\in D^b_c(M^{\K}(\G,\X),\Q_\ell)$,
la flèche composée
$$\begin{array}{rcl}u:d^*i^*Rj_*L={i'}^*\overline{T}_g^*Rj_*L & \stackrel{CB}{\fl} & {i'}^*Rj'_*T_g^*L={i'}^*Rj'_*T_g^!L \\
& \stackrel{CB}{\fl} & {i'}^*\overline{T}_g^!Rj_*L \\
& \fl & d^!i^*Rj_*L=d^*i^*Rj_*L\end{array}$$
(où la dernière flèche est celle définie avant le lemme \ref{lemme:restr_corr_bord})
est égale à $[\Hr_{\lin,r}:\Hr'_{\lin,r}].id$.

\end{proposition}

\begin{corollaire}$r$ est comme dans la proposition ci-dessus,
et on se donne $h_1,h_2,h'\in\G(\widehat{\Z})$, $g_1,g_2\in\G(\Af^p)$
tels que $g_1^{-1}\K'g_1\subset\K$, $g_2^{-1}\K' g_2\subset\K$ et 
$\overline{T}_{g_1}$ (resp. $\overline{T}_{g_2}$)
envoie $Im(i'_{h',r})$ dans $Im(i_{h_1,r})$ (resp. $Im(i_{h_2,r})$).
Il existe donc $q_1,q_2\in\Pa_r(\Q)\QP_r(\Af)\cap\G(\Af^p)$ et $k_1,k_2\in\K$
tels que $q_1h_1=h'g_1k_1$ et $q_2h_2=h'g_2k_2$,
et les diagrammes suivants sont commutatifs
$$\xymatrix{M^{\K'_r}(\G_r,\X_r)\ar[r]^{i'_{h',r}}\ar[d]_{T_{\overline{q}_1}} & M^{\K'}(\G,\X)^*\ar[d]_{\overline{T}_{g_1}} & M^{\K'_r}(\G_r,\X_r)\ar[r]^{i'_{h',r}}\ar[d]_{T_{\overline{q}_2}} & M^{\K'}(\G,\X)^*\ar[d]_{\overline{T}_{g_2}} \\
M^{\K_r}(\G_r,\X_r)\ar[r]^{i_{h_1,r}} & M^{\K}(\G,\X)^* & M^{\K_r}(\G_r,\X_r)\ar[r]^{i_{h_2,r}} & M^{\K}(\G,\X)^*}$$
où $\overline{q}_1$ et $\overline{q}_2$ sont les images dans $\G_r(\Af^p)$
de $q_1$ et $q_2$.

Soit $V\in D^b(Rep_{\G(\Q_\ell)})$.
Alors la correspondance cohomologique
$v_1:T_{\overline{q}_1}^*i_{h_1,r}^*Rj_*\F^{\K}V\fl T_{\overline{q}_2}^*i_{h_2,r}^*Rj_*\F^{\K}V$
définie par
$$\begin{array}{rcl}T_{\overline{q}_1}^*i_{h_1,r}^*Rj_*\F^{\K}V & = & {i'_{h',r}}^*\overline{T}_g^*Rj_*\F^{\K}V \\
& \stackrel{Rj_*u_g}{\fl} & {i'_{h',r}}^*\overline{T}_1^!Rj_*\F^{\K}V \\
& \fl & T_{\overline{q}_2}^!i_{h_2,r}^*Rj_*\F^{\K}(V) \\
& = & T_{\overline{q}_2}^*i_{h_2,r}^*Rj_*\F^{\K}V\end{array}$$
est égale à $[\Hr_{\ell,r}:\Hr_{\ell,r}']$ fois la correspondance
\begin{flushleft}$\displaystyle{v_2:T_{\overline{q}_1}^*i_{h_1,r}^*Rj_*\F^{\K}V\simeq
\F^{\K'_r}R\Gamma(\Gamma_{\ell,r},R\Gamma(Lie(\N_r),q_1h_1.V))}$\end{flushleft}
\begin{flushright}$\displaystyle{\stackrel{q_2^{-1}h_2^{-1}q_1h_1}
{\overrightarrow{\hspace{2cm}}}
\F^{\K'_r}R\Gamma(\Gamma_{\ell,r},R\Gamma(Lie(\N_r),q_2h_2.V))\simeq T_{\overline{q}_2}^*i_{h_2,r}^*Rj_*\F^{\K}V.}$\end{flushright}

\end{corollaire}

\begin{proof} Notons
$$\begin{array}{rcl}\varphi_j:T_{\overline{q}_j}^*i_{h_j,r}^*Rj_*\F^{\K}V & = & {i'_{h',r}}^*\overline{T}_{g_j}^*Rj_*\F^{\K}V \\
&\stackrel{CB}{\fl} & {i'_{h',r}}^*Rj'_*T_{g_j}^*\F^{\K}V \\
& \simeq & {i'_{h',r}}^*Rj'_*\F^{\K'}g_j.V \\
& \stackrel{g_j}{\fl} & {i'_{h',r}}^*Rj'_*\F^{\K'}V\end{array}$$
pour $j=1,2$, et
$$\begin{array}{rcl}\psi:{i'_{h',r}}^*Rj'_*\F^{\K'}V & \stackrel{g_2^{-1}}{\fl} & {i'_{h',r}}^*Rj'_*\F^{\K'}g_2.V \\
& \simeq & {i'_{h',r}}^*Rj'_*T_{g_2}^!\F^{\K}V \\
& \stackrel{CB}{\fl} & {i'_{h',r}}^*\overline{T}_{g_2}^!Rj_*\F^{\K}V \\
& \fl & T_{\overline{q}_2}^!i_{h_2,r}^*Rj_*\F^{\K}V\end{array}$$

On a $v_1=\psi\varphi_1$ par définition de $Rj_*u_g$,
$\varphi_2 v_2=\varphi_1$ d'après la proposition \ref{Pink:corr_Hecke}
et $\psi\varphi_2=[\Hr_{\ell,r}:\Hr'_{\ell,r}]id$ d'après la proposition ci-dessus,
donc 
$$v_1=\psi\varphi_1=\psi\varphi_2 v_2=[\Hr_{\ell,r}:\Hr'_{\ell,r}]v_2.$$

\end{proof}

\begin{proof}[Démonstration de la proposition.] 
On peut supposer que $g=1$, et on le fera pour simplifier les notations.

Comme $\overline{T}_1$ est un morphisme fini entre schémas normaux, on dispose
d'un morphisme fonctoriel trace $Tr_{\overline{T}_1}:\overline{T}_{1*}\overline{T}_1^*\fl id$
(SGA XVII 6.2.5 et 6.2.6).
On notera encore $Tr_{\overline{T}_1}$ le morphisme fonctoriel
$\overline{T}_1^*\fl\overline{T}_1^!$ qui s'en déduit par adjonction.

D'après le lemme \ref{lemme:unicite_prolong}, le morphisme
$$\overline{T}_1^*Rj_*L\stackrel{CB}{\fl}Rj'_*T_1^*L=Rj'_*T_1^!L\stackrel{CB}{\fl}\overline{T}_1^!Rj_*L$$
est égal à
$$Tr_{\overline{T}_1}:\overline{T}_1^*Rj_*L\fl\overline{T}_1^!Rj_*L.$$

Soit $Z=(M_1\times_{M^{\K}(\G,\X)^*}M^{\K'}(\G,\X)^*)_{red}$
(c'est l'union disjointe des strates $Im(i'_{b,r})$ de $M^{\K'}(\G,\X)^*$
au-dessus de $M_1$).
On a un diagramme commutatif dont le carré du bas est cartésien aux nilpotents près
$$\xymatrix@C=15pt{M_1'\ar[rrrd]^{i'}\ar[rd]^k\ar[ddr]_d & & & \\
 & Z\ar[rr]^-I\ar[d]^c & & M^{\K'}(\G,\X)^*\ar[d]^{\overline{T}_1} \\
 & M_1\ar[rr]_-i & & M^{\K}(\G,\X)^*}$$
où $c$ est fini étale et $k$ est une immersion ouverte.

Notons $n_T$ la pondération de $\overline{T}_1$ définie dans SGA XVII 6.2.6 :
si $x'$ est un point géométrique de $M^{\K'}(\G,\X)^*$ d'image $x$ dans $M^{\K}(\G,\X)^*$,
et si $K_{x'}$ (resp. $K_{x}$) est le corps des fractions de l'anneau strictement local
de $M^{\K'}(\G,\X)^*$ en $x'$ (resp. de $M^{\K}(\G,\X)^*$ en $x$),
alors $n_T(x')=[K_{x'}:K_x]$.
En particulier, si l'image de $x'$ est dans $M^{\K'}(\G,\X)$, alors $n_T(x')=1$
(puisque $T_1$ est étale).
Pour tout point géométrique $x$ de $M^{\K}(\G,\X)$, on a donc 
\begin{equation*}
\tag{+}
{\sum_{x'\in\overline{T}_1^{-1}(x)}n_T(x')=deg(T_1)=[\K:\K'].}
\end{equation*}

$n_T$ donne par changement de base une pondération de $c$, d'où des morphismes trace
$Tr_{c}:c_*c^*\fl id$ et $Tr_c:c^*\fl c^!$.
D'après la compatibilité aux changements de base du morphisme trace, 
le morphisme fonctoriel
$$c^*i^*=I^*\overline{T}_1^*\stackrel{I^*Tr_{\overline{T}_1}}{\overrightarrow{\hspace{1.5cm}}}I^*\overline{T}_1^!\fl c^!i^*$$
est égal au morphisme $Tr_c:c^*i^*\fl c^!i^*$,
donc l'endomorphisme $u$ de $d^*i^*Rj_*L$ qu'on veut calculer est égal au morphisme
$$d^*i^*Rj_*L=k^*c^*i^*Rj_*L\stackrel{Tr_c}{\fl}k^*c^!i^*Rj_*L=d^*i^*Rj_*L.$$
Pour conclure, il suffit donc de montrer que si $x'$ est un point géométrique
de $M^{\K'}(\G,\X)^*$ qui se factorise par $i':M_1'\fl M^{\K'}(\G,\X)^*$, alors 
$n_T(x')=[\Hr_{\ell,r}:\Hr'_{\ell,r}]$.

Soit $x'$ un tel point. Son image $x$ dans $M^{\K}(\G,\X)^*$ se factorise évidemment par
$i:M_1\fl M^{\K}(\G,\X)^*$. 
D'après la propriété SGA XVII 6.2.4 (*) des pondérations et l'égalité (+) ci-dessus, on a 
$$\sum_{x''\in\overline{T}_1^{-1}(x)}n_T(x'')=[\K:\K'].$$
Or le groupe $\K/\K'$ agit transitivement sur $\overline{T}_1^{-1}(x)$, donc
les $n_T(x'')$, $x''\in\overline{T}^{-1}(x)$, sont tous égaux.
On en déduit que $n_T(x')=[\K:\K']/\card(\overline{T}_1^{-1}(x))$.

Il reste à calculer $\card(\overline{T}_1^{-1}(x))$.
Notons $N$ le nombre de strates de $M^{\K'}(\G,\X)^*$ de la forme $Im(i'_{b,r})$
qui sont envoyées sur $M_1$ par $\overline{T}_1$. 
Comme on a $Im(i'_{b_1,r})=Im(i'_{b_2,r})$ si et seulement si 
$\Pa_r(\Q)\QP_r(\Af)b_1\K'=\Pa_r(\Q)\QP_r(\Af)b_2\K'$, $N$ est égal au cardinal
de la fibre en $\Pa_r(\Q)\QP_r(\Af)h\K$ de l'application évidente\newline
$\protect{\Pa_r(\Q)\QP_r(\Af)\sous\G(\Af)/\K'\fl\Pa_r(\Q)\QP_r(\Af)\sous\G(\Af)/\K}$,
c'est-à-dire à $[\K:\K']/[\Hr_r:\Hr'_r]$.
Comme chaque strate $Im(i'_{b,r})$ de $M^{\K'}(\G,\X)^*$ qui s'envoie sur $M_1$
a $[\K_r:\K'_r]$ points géométriques au-dessus de $x$, on trouve finalement
$$\card(\overline{T}_1^{-1}(x))=[\K_r:\K'_r]\frac{[\K:\K']}{[\Hr_r:\Hr'_r]}=\frac{[\K:\K']}{[\Hr_{\ell,r}:\Hr'_{\ell,r}]},$$
ce qui finit la preuve.

\end{proof}

\end{document}